\documentclass[10pt,a4paper]{article}
\usepackage{bbm}

\usepackage[leqno]{amsmath}
\usepackage{amsfonts}
\usepackage{graphicx}
\usepackage{amsmath}
\usepackage{amssymb}
\usepackage{latexsym}
\usepackage{amsmath, amsfonts,amssymb, amsthm, euscript,makeidx,color,mathrsfs}
\usepackage{enumerate}

\usepackage[colorlinks,linkcolor=blue,anchorcolor=green,citecolor=red]{hyperref}%when print, use the next package
\def\sqr#1#2{{\vcenter{\vbox{\hrule height.#2pt
              \hbox{\vrule width.#2pt height#1pt \kern#1pt \vrule width.#2pt}
              \hrule height.#2pt}}}}
\def\5n{\negthinspace \negthinspace \negthinspace \negthinspace \negthinspace }
\def\4n{\negthinspace \negthinspace \negthinspace \negthinspace }
\def\3n{\negthinspace \negthinspace \negthinspace }
\def\2n{\negthinspace \negthinspace }
\def\1n{\negthinspace }

\def\dbE{\mathbb{E}}
\def\dbF{\mathbb{F}}

\def\dbH{\mathbb{H}}

\def\dbN{\mathbb{N}}

\def\dbP{\mathbb{P}}

\def\dbR{\mathbb{R}}

\def\sK{\mathscr{K}}
\def\sL{\mathscr{L}}

\def\sU{\mathscr{U}}

\def\fU{\mathfrak{U}}

\def\={\buildrel \triangle \over =}

\def\ds{\displaystyle}

\def\ns{\noalign{\ss}}
\def\a{\alpha}
\def\b{\beta}
\def\g{\gamma}
\def\d{\delta}
\def\e{\varepsilon}
\def\z{\zeta}

\def\l{\lambda}
\def\m{\mu}

\def\si{\sigma}

\def\f{\varphi}
\def\th{\theta}
\def\o{\omega}

\def\i{\infty}
\def\G{\Gamma}
\def\D{\Delta}
\def\Th{\Theta}

\def\O{\Omega}
\def\cA{{\cal A}}
\def\cB{{\cal B}}

\def\cD{{\cal D}}

\def\cF{{\cal F}}

\def\cH{{\cal H}}
\def\cI{{\cal I}}

\def\cL{{\cal L}}
\def\cM{{\cal M}}

\def\cO{{\cal O}}

\def\cS{{\cal S}}

\def\cU{{\cal U}}

\def\no{\noindent}

\def\ss{\smallskip}
\def\ms{\medskip}

\def\q{\quad}
\def\qq{\qquad}
\def\hb{\hbox}

\def\limsup{\mathop{\overline{\rm lim}}}

\def\lt{\left}
\def\rt{\right}
\def\lan{\langle}
\def\ran{\rangle}
\def\llan{\left\langle}
\def\rran{\right\rangle}

\def\tb{\textcolor{blue}}
\def\rf{\eqref}
\def\esssup{\mathop{\rm esssup}}

\def\wt{\widetilde}
\def\ti{\tilde}
\def\cd{\cdot}
\def\cds{\cdots}

\def\ae{\hbox{\rm a.e.}}
\def\as{\hbox{\rm a.s.}}

\def\les{\leqslant}
\def\ges{\geqslant}

\def\({\Big (}
\def\){\Big )}
\def\[{\Big[}
\def\]{\Big]}

\def\bde{\begin{definition}\label}
\def\ede{\end{definition}}
\def\be{\begin{equation}}
\def\bel{\begin{equation}\label}
\def\ee{\end{equation}}
\def\bt{\begin{theorem}\label}
\def\et{\end{theorem}}
\def\bc{\begin{corollary}\label}
\def\ec{\end{corollary}}
\def\bl{\begin{lemma}\label}
\def\el{\end{lemma}}
\def\bp{\begin{proposition}\label}
\def\ep{\end{proposition}}
\def\bas{\begin{assumption}\label}
\def\eas{\end{assumption}}
\def\br{\begin{remark}\label}
\def\er{\end{remark}}
\def\bex{\begin{example}\label}
\def\ex{\end{example}}
\def\ba{\begin{array}}
\def\ea{\end{array}}
\def\ed{\end{document}}
\def\pf{\begin{proof}}
\def\ef{\end{proof}}
\def\eps{\epsilon}
\def\square#1{\vbox{\hrule\hbox{\vrule height#1%
     \kern#1\vrule}\hrule}}
\def\rectangle#1#2{\vbox{\hrule\hbox{\vrule height#1%
     \kern#2\vrule}\hrule}}

\def\T{[0,T]}

\font\tenbb=msbm10 \font\sevenbb=msbm7 \font\fivebb=msbm5

\newfam\bbfam
\scriptscriptfont\bbfam=\fivebb \textfont\bbfam=\tenbb
\scriptfont\bbfam=\sevenbb

\newtheorem{theorem}{\hskip 1.3em Theorem}[section]
\newtheorem{definition}[theorem]{\hskip 1.3em Definition}
\newtheorem{proposition}[theorem]{\hskip 1.3em Proposition}
\newtheorem{corollary}[theorem]{\hskip 1.3em Corollary}
\newtheorem{lemma}[theorem]{\hskip 1.3em Lemma}
\newtheorem{remark}[theorem]{\hskip 1.3em Remark}
\newtheorem{example}[theorem]{\hskip 1.3em Example}

\makeatletter
   
   \@addtoreset{equation}{section}
\makeatother

\newcommand{\norm}[1]{\left\Vert#1\right\Vert}
\def\inte{\int_0^t}

\usepackage{verbatim}

\begin{document}

\title{Singular backward stochastic Volterra integral equations
in infinite dimensional spaces}
\date{ }

\author{
%Yushi Hamaguchi, \footnote{Graduate School of Engineering Science, Department of Systems Innovation, Osaka University. 1-3, Machikaneyama, Toy-onaka, Osaka, Japan. Email: {\tt hmgch2950@gmail.com}}~~~ %This author is supported by National Natural Science Foundation of China (No. 12171471).}
Tianxiao Wang \footnote{School of Mathematics, Sichuan University, Chengdu, China. Email: {\tt wtxiao2014@scu.edu.cn}} ~~ and ~~~
%Jiongmin Yong \footnote{Department of Mathematics, University of Central Florida, Orlando, USA. Email: {\tt jiongmin.yong@ucf.edu}. This author is supported in part by NSF Grant DMS-1812921.}
Mengliang Zheng \footnote{School of Mathematics, Sichuan University, Chengdu, China. Email: {\tt ml\_zheng@stu.scu.edu.cn}}
}

\maketitle
%\begin{}

%\begin{abstract}
% In this paper, a systematic investigation is carried out for the general solvability of multi-dimensional backward stochastic Volterra integral equations (BSVIEs) with the generators being super-linear in the adjustment variable $Z$. Two major situations are discussed: (i) When the free term is bounded with the dependence of the generator on $Z$ being of ``diagonally strictly'' quadratic growth and being sub-quadratically coupled with off-diagonal components; (ii) When the free term is unbounded having exponential moments of arbitrary order with the dependence of the generator on $Z$ being diagonally no more than quadratic and being {\color{red} independent of off-diagonal components}. Besides, for the case that the generator {\color{red} is} super-quadratic in $Z$, some negative results are presented.

%\end{abstract}

\ms

%\bf
\begin{abstract}
In this paper, the notion of singular backward stochastic Volterra integral equations (singular BSVIEs for short) in infinite dimensional space is introduced, and the corresponding well-posedness is carefully established. A class of singularity conditions are proposed, which not only cover that of fractional kernel, Volterra Heston model kernel, completely monotone kernels, to mention a few, but also happen to be used in the forward stochastic Volterra integral with new conclusions arising.
Motivated by mathematical physics problem such as the viscoelasticity/thermoviscoelasticity of materials, heat conduction in materials with memory, optimal control problems of abstract stochastic Volterra integral equations (including fractional stochastic evolution equations and stochastic evolutionary integral equations) are presented. At last, our BSVIEs are surprisingly used in maximum principle of controlled stochastic delay evolution equations. One advantage of this new standpoint is that the final cost functional can naturally depend on the past state for the first time. 	
\end{abstract}

\ms


\section{Introduction}
Throughout this paper, $T>0$ is fixed, $H, V$ are two separable Hilbert spaces and $\cL_2^0 \= \cL_2(V;H).$
Let $(\O,\cF,\dbF \= \{\cF_t\}_{t\in \T},\dbP)$ be a complete filtered probability space on which a  $V$-valued cylindrical Brownian motion $W(\cd)$ is defined with $\dbF$ being its natural filtration augmented by all the $\dbP$-null sets in $\cF$.

The main topic in this paper is to study the following infinite dimensional backward stochastic Volterra integral equation (BSVIE, for short),
%\
\bel{Type-II}
Y(t) = \psi(t) + \int_t^T g(t,s,Y(s),Z(t,s),Z(s,t))ds - \int_t^T Z(t,s)dW(s), \qq t\in \T,
\ee
where $\psi(\cd)$ and $g(\cd)$ are called the  {\it free term} and {\it generator} of \rf{Type-II}, respectively. They are given maps, valued in $H$ and satisfy some singular assumptions specified later. What we are interested in is a pair of adapted process $(Y(\cd),Z(\cd,\cd))$ which satisfies \rf{Type-II} in the usual It\^o's sense and some additional conditions (see Section 3 for details).

When $\psi, g, Z$ are independent of $t$, the
above BSVIE \rf{Type-II} is reduced to the following well-known backward stochastic differential equation (BSDE for short):
\bel{BSDEs}\ba{ll}
\ns\ds Y(t)=\xi+\int_t^T g(s,Y(s),Z(s)) ds-\int_t^T Z(s)dW(s),
\ea\ee
which is introduced by Pardoux-Peng \cite{PP} and extensively studied by many researchers.
Another special case of \rf{Type-II} with $g(\cd)$ independent of $Z(s,t$) and $\psi(\cd) \equiv \xi \in \dbR^n$ was firstly studied by Lin \cite{Lin}, and followed by e.g. Aman-N'Zi \cite{AN}, Wang-Zhang \cite{Wang-Zhang 2007}, Djordjevi\'c-Jankovi\'c \cite{DJ}, Hu-Oksendal \cite{HuO}, etc.
BSVIEs of \rf{Type-II} are inspired by optimal control problems for forward SVIEs and were first studied by Yong \cite{Yong 06}. Based on the previous work \cite{Yong 06},  Yong \cite{Y 2008} continued the study of \rf{Type-II} in a more systematic way with certain singular Lipschitz functions and proposed a new concept of solution: adapted M-solution.
After the work of Yong \cite{Yong 06, Y 2008}, BSVIEs have attracted
many researchers' interest and are extended in various forms, such as backward doubly SVIEs (Shi et al. \cite{SWX}), 
%\tb{mean-field backward doubly SVIEs (Wu-Hu \cite{WHu}),}
reflected BSVIEs (Agram-Djehiche \cite{Agram-Djehiche}), time delayed BSVIEs (Bessner-Rosazza Gianin \cite{BRG}), BSVIEs with diagonal-solution generator (Hernandez-Possamai \cite{Hernandez-Possamai 2021}, Hernandez \cite{Hernandez 2021}, Wang-Yong \cite{WHY}), backward stochastic Volterra integro-differential equations (Wang \cite{Wang-2022}), mean-field BSVIEs (Shi et al \cite{Shi-Wang-Yong 2013}), path-dependent BSVIEs (Overbeck and R\"oder \cite{OR}), infinite horizon BSVIEs (Hamaguchi \cite{Hamaguchi1}), BSVIEs with jumps in general filtration (Popier \cite{Popier 2021}). Besides the theoretical aspects, BSVIEs also have some interesting applications.
\begin{itemize}
\item Stochastic control theory: As is mentioned before, backward stochastic Volterra integral equations of form \rf{Type-II} are originally motivated by the study of maximum principle for optimal control of SVIEs. After the seminar work \cite{Yong 06}, there are several following-up works,  Agram-Oksendal \cite{AO}, Shi-Wang-Yong \cite{SWY}, and Wang \cite{W}, Wang-Yong \cite{Wang-Yong-2023}, Wang-Yong-Zhou \cite{WHYZhou}, Wang-Zhang \cite{WZ}, to mention a few.

\item Mathematical finance:
It is well known that many financial problems can be formulated by BSDEs, such as recursive utility, dynamic risk measures, asset pricing in incomplete markets, and so on.
Similarly, it is found similar topics can also be discussed via the BSVIEs, such as Yong \cite{Yong 2007}, Wang-Sun-Yong \cite{WSY}, and Beissner--Rosazza Gianin \cite{BRG}.

\item PDEs theory: BSVIEs have close connections with partial differential equations. In Wang-Yong \cite{WY}, they established a representation of adapted M-solutions to BSVIEs in terms of the solution to a system of (non-local) partial differential equations. Following-up works along this line can be found in e.g. Wang \cite{WH}, Wang-Yong-Zhang \cite{WHYZhang}, Lei-Pun \cite{LP}.

\item Time-inconsistent control problems: In real world, people's subjective time-preferences and risk-preferences usually lead to the time-inconsistency phenomenons. Fortunately, BSVIEs can well capture these characters and are used to study the time-inconsistent stochastic optimal control problems, see e.g. Wang-Yong \cite{WHY}, Wang-Zheng \cite{Wang-Zheng-2021} or the survey of Yan-Yong \cite{YY}. We emphasize that this area indicates an \it interesting application advantage of BSVIEs in contrast to the BSDEs. \rm
\end{itemize}

\ms

The aim of the current paper is to introduce and establish the \it singular \rm backward stochastic Volterra integral equations \rm theory in the \it infinite dimensional \rm framework. To our best knowledge, this notion seems to be new in the literature. In contrast with the existing literature on BSVIEs, we highlight two keywords in the current study, i.e., singularity and infinite dimensional setting.

The motivation to study BSVIEs in the singular scenario is twofold. On the one hand,
Volterra integral equations have close connections with singularities from their beginning. In fact, the earliest research can be traced back to 1823 when Abel proposed a generalization of the tautochrone problem formulated as a linear singular convolution integral equation (\cite{Abel}). Since then, singular (forward) Volterra integral equations have attracted a lot of attention, such as the explicit expression of the solution \cite{D}, the asymptotic solution
\cite{HO}, the smoothness of the solution \cite{MF}, the numerical solution  \cite{L}, and so on. Therefore, our study can be seen as a continuation in both the \it backward, stochastic \rm setting.
On the other hand, from an application point of view, the singularity in the Volterra integral equations can be used to describe some memory or hereditary phenomena arising from mathematic physics (\cite{HO, RM}) or mathematical finance (\cite{BD, EER}). For example, it has been proved that singular Volterra integral equations can fit remarkably well in the so-called rough volatility models in financial markets, see \cite{BD, EER}. In addition, it was found that (time) fractional differential equations are closely related to singular Volterra integral equations (e.g. \cite{LY}). If we combine the time fractional derivative with the well-known BSDEs \rf{BSDEs}, we then arrive at the following example, which is seen as a particular form of singular BSVIE.
\bex{NA} Consider the following Caputo fractional backward stochastic differential equations of order $\a \in (\frac{1}{2} , 1)$ on $[0, T ]$ (see e.g. \cite{NA}):
\bel{CfBSDE}
\left\{\begin{aligned}
	&\left( _{t}^{C} D^{\alpha}_T y\right)(t)=A y(t)-f(t, y(t), z(t, s))-z(t, s) \frac{d w(t)}{d t}, \q t\in [0,T],\\
	&y(T)=\xi,
\end{aligned}\right.
\ee
\no where $^C_tD^\a_T$ denotes the right Caputo fractional derivative of order $\a$ defined by
$$_t^CD^\a_T y(t)\= -\frac{1}{\G(1-\a)} \int_t^T \frac{y'(s)}{(s-t)^\a}ds, $$
$A$ is a constant matrix and $(w(\cd))_{t\ges 0}$ is a standard $m$-dimensional
Brownian motion defined on a complete probability space $(\O,\cF,\dbF,\dbP)$ with natural filtration $\dbF$ satisfying usual conditions, $\frac{dw}{dt}$ denotes the white noise, which is the generalized derivative of Browian motion. Moreover,  $\xi \in L_{\cF_T}^2(\O,\dbR^n)$, $f:[0,T] \times \dbR^n \times \dbR^{n\times m} \to \dbR^n $ is measurable such that
$f(\cd,0,0)\in L^2(0,T;\dbR^n)$ and for all $y, \bar y\in \dbR^n, z,\bar z \in \dbR^{n\times m}, t\in [0,T]$, there exists $L>0$ satisfying
$$ |f(t,y,z)- f(t,\bar y,\bar z)|^2 \les L\( |y-\bar y|^2 + |z-\bar z|^2 \), \q a.s.$$
According to \cite{NA}, a pair of $(x, y)\in  M[0, T]\= L_\dbF^2(0,T;\dbR^n) \times L_\dbF^2(\D^*; \dbR^{n\times m})$ is named as a mild solution of \rf{CfBSDE} if it satisfies
\bel{mfbsde} \begin{aligned}
	y(t) = &\ \xi + \frac{1}{\G(\a)} \int_t^T (s-t)^{\a-1} \[f(s,y(s),z(t,s)) - Ay(s)\] ds \\
	&+ \frac{1}{\G(\a)} \int_t^T (s-t)^{\a-1} z(t,s) dw(s), \q t\in [0,T].
\end{aligned}\ee
If we denote
$$Y(t)\=y(t), \q Z(t,s)\=(s-t)^{\a-1} z(t,s),$$
then \rf{mfbsde} can be seen as a special case of BSVIE \rf{Type-II} with
$$\left\{ \begin{aligned}
	&\psi(t) \= \xi, \\
    &g(t,s,Y(s),Z(t,s)) \= (s-t)^{\a-1} \[ f(s,Y(s),(s-t)^{1-\a}Z(t,s)) - AY(s) \].
\end{aligned}\rt.$$
%
%Since $A$ can be viewed as a bounded linear operator, \rf{CfBSDE} then admits a unique mild solution according to Theorem \ref{43}. %{\color{red} Check the condition in (H3.1)?}
\ex

After the discussion of singularity, let us look at the other keyword in the current study: infinite dimensional setting. We point out that our research motivation comes from mathematical physics. To this end, we present the following example of BSVIEs \rf{Type-II} relevant to the stochastic evolutionary integral equations that describe the viscoelasticity, thermoviscoelasticity of material, the incompressible fluids, as well as heat conduction problem in materials with memory.
\bex{}
Inspired by e.g. \cite{P, Zhang}, we first make some introductions on the forward semilinear stochastic evolutionary integral equation.
Suppose $H$ is a Hilbert space, $A$ is a closed linear unbounded operator on $H$ with dense domain $\cD(A)$, and $a\in L^1(0,T;\dbR_+)$ is a scalar kernel. Consider
\bel{seie}
X(t)=x_0-\int_0^ta(t-s)AX(s)ds+\int_0^t\Phi\big(s,X(s)\big)ds+\int_0^t\Psi\big(s,X(s)\big)dW(s),\ \ t\in[0,T],
\ee
where $\Phi: [0,T] \times \O \times H \to H$ and $\Psi: [0,T] \times \O \times H \to \cL_2^0$ are Lipschitz and linear growth. As is shown later, with different form of $a(\cd)$, such kind of equations have interesting applications in mathematical physics.

To study the well-posedness of (\ref{seie}), we need the notion of
resolvent (\cite{P}). A family of bounded linear operators in $H$, denoted by $\{S(t)\}_{t\ges0}$, is called a resolvent for \rf{seie} if $t\mapsto S(t)$ is strong continuous and commutes with $A$, $S(0)=I$, and the following resolvent equation holds:
$$
S(t)x = x- \int_0^t a(t-s) A S(s)x ds,\ \ x\in\cD(A),\ \ t\in[0,T].
$$
We point out several cases to ensure the existence of resolvent, as well as applications in various frameworks.

\begin{itemize}
		\item  If $a(t)=1$ and $a(t)=t$, the above resolvent $S(\cd)$ becomes the $C_0$-semigroup generated by $A$, and the cosine family $\hbox{Co}(t)$ generated by $A$, resp.

		\item  If $a(\cd)$ is a Bernstein function, $A$ generates a bounded cosine family $\hbox{Co}(\cd)$, then the resolvent exists and is bounded.  It is applicable to Stokes' first problem (or the Rayleigh problem) in viscoelasticity and torsion problem of a rod.

\item If $a(\cd)$ is $k$-regular function (see \cite[Definition 3.3]{P}), $A$ generates a proper analytic $C_0$-semigroup, then the bounded resolvent exists. It is applicable to heat conduction problem in materials with memory.

\item If $a(\cd)$ is completely positive (see e.g. \cite[Definition 4.5]{P}), and $A$ generates a proper analytic $C_0$-semigroup, then the bounded resolvent $S(\cd)$ exists. It is  applicable to models describing isotropic incompressible fluids that cover the linear Navier-Stokes system.

\item If $a$ is a proper creep function (see e.g. \cite[Definition 4.4]{P}), and $A$ is proper selfadjoint and negative, then the bounded resolvent exists. It is applicable to models describing three dimensional isotropic synchronous materials.

\item If $a(\cd)$ is completely monotonic, and $A$ generates a $C_0$-semigroup of contractions (such as classical thermoelasticity operator), then a bounded resolvent S(t) exists. It is applicable to synchronous isotropic thermoviscoelasticity.

	\end{itemize}
We refer to \cite[Theorem 3.1, Theorem 3.2, Theorem 4.2, Theorem 4.3, Theorem 4.4, Chapter 5]{P} for more relevant details on the above resolvent.
Once the bounded resolvent is given, by a (mild) solution of \rf{seie} we mean that $X(\cd)$ satisfies the following stochastic Volterra integral equation:
\bel{Stochastic-evolution-integral-1}\ba{ll}
\ns\ds
X(t) = S(t) x_0 + \int_0^t S(t-s) \Phi\( s,X(s)\) ds + \int_0^t S (t-s) \Psi\( s,X(s)\) dW(s), \q t\in [0,T].
\ea\ee
It can be viewed as a special case of \rf{fSVIE} with
$$\f(t)\= S(t) x_0, \q
A(t,s,x)\= S(t-s) \Phi( s,x),\q
B(t,s,x)\= S(t-s) \Psi( s,x).$$
%
%
% e

As indicated in \cite{Y 2008}, when studying the optimal control problems for \rf{Stochastic-evolution-integral-1}, the following kind of linear BSVIE is required:
\bel{ad2}
Y(t) = \f(t)+ \int_t^T [ M_1(t)^* S (s-t)^* Y(s)+ M_2(t)^*S (s-t)^*Z(s,t) ] ds - \int_t^T Z(t,s) dW(s), \q t\in [0,T],
\ee
with proper operator-valued functions $M_1$, $M_2$. Obviously, the equation \rf{ad2} can be viewed as a special case of \rf{Type-II}.
\ex

From the above arguments, we see that it is an important and interesting topic to consider singular BSVIEs in infinite dimensional spaces from both theoretical and applied standpoints.
We admit that similar singularity thoughts appeared in for example \cite{Y 2008} with finite dimensional case, and  Anh et al. \cite{Y 2011}, Ren \cite{Ren} with infinite dimensional case. In contrast, our singularity framework imposed on  \rf{Type-II} is inspired by Hamaguchi-Wang \cite{Hamaguchi}. As is shown in Remark \ref{rmk41}, our singularity assumptions, even for the strongest one, are weaker than theirs. Meanwhile, it is found that similar singularities can be well fit to the forward stochastic Volterra system, which help us extend and improve the existing singularity conditions in \cite{Abi1,Abi2}, see also the Remark \ref{comparison-SVIE}. We also make detailed comparisons with two important papers \cite{LY, Zhang} in the Appendix.
As applications, we study the optimal control problem for singular forward stochastic Volterra integral equations (including fractional  stochastic evolution equations and controlled stochastic evolutionary integral equations as special cases)
and use linear singular BSVIEs as the adjoint equation to derive the Pontragin maximum principle. %
Eventually, we give another interesting and useful applications to optimal control problem of stochastic delay evolution equation where the terminal cost is allowed to depend on the past of the state.

\medskip

At this moment, let us list the main novelties of the current paper as follows:
\begin{itemize}
\item  Singularity: In contrast with the existing literature on BSVIEs, we give some new and deeper treatment on the singularity of the generator. First, we carefully impose different singular assumptions on each element of $g$ in \rf{Type-II} (see (H3.1) and Remark \ref{remark37}) to highlight their differences. Second,  the corresponding singular kernels we propose is non-convolution which enable us to cover other Rough, fractional kernels, (see Remark \ref{comparison-SVIE}). Third, we apply such singularities idea to the forward systems and obtain their well-posedness that seems to be new to our best knowledge.

\item Volterra system feature: To show the advantage of Volterra system, we present one optimal control problem of controlled stochastic delay evolution equations. More precisely, to obtain the maximum principle of optimal controls, we rewrite the variational equation as a forward SVIEs and thus use linear BSVIEs as the adjoint equation instead of anticipated BSDEs. One advantage is that the final cost can easily and naturally depend on the past of the state, which in some sense overcome the limitation of anticipated BSDEs in such general scenario.

\item Infinite dimensionality: The infinite dimensional feature of equations \rf{Type-II} enables us to apply the current consideration to other various stochastic PDEs arising from mathematical physics (say, the issues of viscoelasticity, thermoviscoelasticity, the incompressible fluids, heat conduction in materials with memory), such as the fractional stochastic evolution equations, stochastic delay evolution equation, stochastic evolutionary integral equations and so on.

\end{itemize}

The rest of this paper is organized as follows. Section 2 collects some preliminary results.  We give the well-posedness of singular BSVIEs in section 3. %The existence and uniqueness of adapted solution and adapted M-solution for Type-I and Type-II BSVIEs are investigated in section 4 and section 5, respectively.
As applications, we study two classes of optimal control problems with convex control regions in Section 4, respectively. In Section 5, we give a concluding remark. And at last, some comparison results are presented in the Appendix.

%{\color{blue}arise frequently in the continuum mechanics for materials with memory, i.e. the theory of viscoelastic materials, for example, some standard problems, the Rayleigh problem and torsion of a rod. Besides, the heat conduction in materials with memory and the viscoelastic beams and plates can be lead to evolutionary integral equations as shown in \cite{P}.}

\section{Preliminaries}
Let $V, H$ be two separable Hilbert spaces. {Denote by $\cL(V;H)$ the space of all bounded linear operators from $V$ to $H$} and  $\cL_2^0 $ the space of all Hilbert-Schmidt operators from $V$ to $H$, i.e.
$$\cL_2^0 \= \lt\{ F \in \cL(V;H) \;\Big| \sum_{i=1}^\i |Fe_i|_H^2<\i, \text{ where } \{e_i\}_{i=1}^\i \text{ is an orthonormal
	basis of } V\rt \}. $$
One can show that, $\cL_2^0$ equipped with the inner product
$$ \lan F,G\ran_{\cL_2^0} \= \sum_{i=1}^\i \lan Fe_i, Ge_i \ran_H, \q \forall F, G \in \cL_2^0,$$
is a separable Hilbert space.

Next, we define the triangle domain
$$\D\=\big\{(t,s)\in[0,T]^2\bigm|0\les s<t \les T\big\},$$
and for each $R\in[0,T)$,
$$\D[R,T]\=\big\{(t,s)\in[R,T]^2\bigm|0\les R \les s<t \les T\big\}.$$
Similarly, we define another triangle domain
$$\D^*\=\big\{(t,s)\in[0,T]^2\bigm|0\les t < s \les T\big\},$$
and for each $R\in[0,T)$,
$$\D^*[R,T]\=\big\{(t,s)\in[R,T]^2\bigm|0\les R\les t < s\les T\big\}.$$

Define the space of $H$-valued square integrable random variables as follows:
$$L^2_{\cF_t}(\Omega;H)\=\Big\{\xi:\Omega\to H\bigm|\xi\hb{ is $\cF_t$-measurable, }\|\xi\|_2\equiv\big(\dbE|\xi|_H^2\big)^{1\over 2}<\infty\Big\}.$$
Obviously,  $L^2_{\cF_t}(\Omega;H)$ is a Banach space under the norm $\|\,\cd\,\|_2$. %When the range space $H$ is clear from the context and is not necessarily emphasized, we will omit $H$. In particular, we will denote $L^2_{\cF_T}(\Omega)=L^2_{\cF_T}(\Omega;H)$.

\ms

Next, we introduce spaces of stochastic processes. In order to avoid repetition, all processes $(t,\o)\mapsto\f(t,\o)$ are at least $\cB([0,T])\otimes\cF_T$-measurable without further mentioning, where $\cB([0,T])$ is the Borel $\si$-field of $[0,T]$. For any $0 \les R \les S \les T$, we define
$$\ba{ll}
\ds L_{\cF_S}^2(R,S;H)\1n\=\1n\Big\{\f:[R,S]\1n\times\1n\O\to H\bigm| \f \text{ is } \cB([R,S])\otimes\cF_S \text{-measurable and }
\dbE\(\int_R^S\3n|\f(t)|_H^2dt\) \3n<\infty\Big\},\\

\ds
L_\dbF^2 (R,S;H) \= \Big\{ \f \in L_{\cF_S}^2(R,S;H) \bigm| \f \text{ is $\dbF$-adapted} \Big\}, \\

\ds
L_\dbF^2(\O;C([R,S];H)) \= \Big\{ \f: [R,S] \times \O \to H \bigm| \f(\cd) \text{ is $\dbF$-adapted, continuous and } \dbE\(\max_{R \les t \les S}|\f(t)|_H^2 \) < \i \Big\},
\\

\ds
L_\dbF^2(R,S;\cL_2^0) \= \Big\{ \f: [R,S] \times \O \to \cL_2^0 \bigm| \f(\cd) \text{ is $\dbF$-adapted and } \dbE\( \int_R^S|\f(\cd)|^2_{\cL_2^0}dt \) < \i \Big\}, \ea$$
\\
%L^2(R,S;L_\dbF^2(S,T;\cL_2^0)) = \\
%\ea$$
%
%$$\ba{ll}
%
and $$\ba{ll}
%\ns\ds L^2_\dbF(\O;L^2(\D^*;H))\=\Big\{\z:\D^*\times\Omega\to H\bigm|\forall r\in[0,T],\, \z(r,\cd)\in L^2_\dbF(\O;L^2(r,T;H)),\\
%
%\ns\ds\qq\qq\qq\qq\qq\qq\qq\qq\qq\qq\qq\dbE\int_0^T\(\int_r^T|\z(r,s)|_H^2ds\)dr<\infty\Big\},\\
%
\ns\ds L^2(R,S;L_\dbF^2(S,T;\cL_2^0))\=\Big\{\z:[R,S]\times [S,T] \times \Omega\to \cL_2^0 \bigm|\forall r\in[R,S],\, \z(r,\cd)\in L_\dbF^2(S,T;\cL_2^0), \text{ and }\\
\ns\ds\qq\qq\qq\qq\qq\qq\qq\qq\qq\qq\qq\dbE\int_R^S\(\int_S^T|\z(r,s)|_{\cL_2^0}^2ds\)dr<\infty\Big\},\\
\ns\ds L^2(R,S;L^2_\dbF(R,S;\cL_2^0))\=\Big\{\z:[R,S]^2\times \Omega\to \cL_2^0 \bigm|\forall r\in[R,S],\, \z(r,\cd)\in L^2_\dbF(R,S;\cL_2^0), \text{ and }\\
\ns\ds\qq\qq\qq\qq\qq\qq\qq\qq\qq\qq\qq\dbE\int_R^S\(\int_R^S|\z(r,s)|_{\cL_2^0}^2ds\)dr<\infty\Big\},\\
\ds
\cH^2[R,S] \= L_\dbF^2(R,S;H) \times L^2(R,S;L_\dbF^2(R,S;\cL_2^0)).
\ea$$

\ms

\no Define the following spaces:

$$L^2(\D^*;\dbR_+) \= \lt\{ f:\D^*\mapsto \dbR_+ \Big|\; f \text{ is  measurable and } \int_0^T \int_t^T |f(t,s)|^2dsdt <\i \rt\},$$
$$L^2(\D;\dbR_+) \= \lt\{ f:\D\mapsto \dbR_+  \Big|\; f \text{ is  measurable and } \int_0^T \int_s^T |f(t,s)|^2dtds   <\i \rt\},$$

\no and $\mathscr{L}^{2}\left(\D^*; \mathbb{R}_+\right)$ the set of $f \in L^{2}\left(\D^*; \mathbb{R}_+\right)$ satisfying the following two conditions: \\
1) it holds that
$$ \|f\|_{\mathscr{L}^{2}\left(\Delta^*;\dbR^+\right)}\=\esssup_{t \in\left(0, T\right)}\left(\int_t^T |f(t,s)|^{2} ds\right)^{1 \over 2}<\infty; $$
2) for any $\e > 0$, there exists a finite partition $\left\{U_{i}\right\}_{i=0}^{m}$ of $\left(0, T\right)$ with $0=U_{0}<U_{1}<\cdots<U_{m}=T$ such that
$$ \esssup_{t \in\left(U_i, U_{i+1}\right)}\left(\int_t^{U_{i+1}} |f(t,s)|^{2} ds\right)^{1 \over 2}<\varepsilon $$
for each $ i \in\{0,1, \ldots, m-1\} $.\\
The space $\sL^2(\D,\dbR_+)$ can be defined similarly. For notation simplicity, in the following we omit the range space $\dbR_+$ without confusion, i.e.
$$  L^2(\D^*) \= L^2(\D^*;\dbR_+),  L^2(\D) \= L^2(\D;\dbR_+),   \sL^{2}(\D^*) \= \sL^{2}\left(\D^*; \mathbb{R}_+\right), \sL^2(\D) \=\sL^2(\D,\dbR_+).$$
\medskip

%At last, we introduce the following stochastic process spaces.
%
The following example shows the relation between the two conditions in $\sL^2(\D^*)$, from which we see that neither 1) nor 2) is redundant and can be dropped. %To sum up, one can not drop any one of them.

\bex{r} (i) Inspired by Shi and Wang \cite{W 2012}, let $\displaystyle f(t,s) \= \sqrt{\frac{2}{T-t}}, \q (t,s)\in \D^*$. It is easy to see that
$$\esssup_{t\in (0,T)} \int_t^T f(t,s)^2ds =2 <\i.$$
Therefore, condition 1) in the definition is fulfilled.
However, for any partition $\left\{U_{i}\right\}_{i=0}^{m}$ of $\left(0, T\right)$ with $0=U_{0}<U_{1}<\cdots<U_{m}=T$, by the choice of $f(\cd,\cd)$,
$$\esssup_{t\in (U_{m-1},T)} \int_t^T f(t,s)^2ds =2,$$
which shows the condition 2) does not hold any more. In other words, condition 1) can not imply condition 2), and hence condition 2) can not be removed.

(ii) On the other hand, let us look at $f(t,s) = f_1(s-t), (t,s) \in \D^*$, with $f_1 \in L^2(0,\frac{T}{2}) \setminus L^2(0,T)$. %for example, let $\ds f(t) = \sqrt{\frac{1}{T-t}}, \, t\in (0,T).$
In this case, for any $\e>0$, there exist $\d>0$ and a partition $\{U_i\}_{i=0}^m$ of $\left(0, T\right)$ with $0=U_{0}<U_{1}<\cdots<U_{m}=T$, whenever $\max\limits_{i\in \{0,1,\cds,m-1\}} (U_{i+1} - U_i) <\min\{\d,\frac{T}{2}\}$,
\be \notag \begin{aligned}
	\esssup_{t \in (U_i,U_{i+1})}	\int_t^{U_{i+1}} |f(t,s)|^2ds &= \esssup_{t \in (U_i,U_{i+1})}\int_0^{U_{i+1} -t}|f_1(s)|^2ds \\
	&= \int_0^{U_{i+1} - U_i}  |f_1(s)|^2ds <\e, \quad \forall i \in \{0,1,\cds,m-1\}.
\end{aligned} \ee
However,
\be \notag \begin{aligned}
	\esssup_{t \in (0,T)}	\int_t^{T} |f(t,s)|^2ds &= \esssup_{t\in (0,T)}\int_0^{T -t}|f_1(s)|^2ds = \int_0^{T}  |f_1(s)|^2ds =\i,
\end{aligned} \ee
which shows that condition 1) does not hold. In other words, condition 2) does not imply condition 1), and hence condition 1) can not be dropped.
\ex
The following gives example in $L^2(\D^*)$ and $\sL^2(\D^*)$ respectively, from which we can clearly see the difference between these two spaces.
\bex{d-s-k}
  For $\a,\b \in [0,1)$, let
  $$k^{\alpha,\beta}(t,s)\=\frac{1}{(s-t)^\alpha t^\beta},\quad (t,s)\in \D^*.$$
  It is clear that when $\a \in [0,\frac{1}{2})$, it holds that
  $$ \int_t^T \(\frac{1}{(s-t)^\alpha t^\beta}\)^2 ds =\frac{1}{1-2\a} \frac{(T-t)^{1-2\a}}{t^{2\b}}.$$
  From this, we have that
  $$\left\{\begin{aligned}
  	 &k^{\a,\b} \in L^2(\D^*), \qq \text{ if } \a \in [0,\frac{1}{2}), \b \in [0,\frac{1}{2}); \\
  	 &k^{\a,\b} \in \sL^2(\D^*), \qq \text{ if } \a \in [0,\frac{1}{2}), \b=0.
  \end{aligned}\right.$$
In the literature, such $k$ is referred to as  a doubly singular kernel.
\ex

Here and next $C$ is a generic positive constant which may be different from line to line.

\ms

\section{Singular backward stochastic Volterra systems}

In this section, we are devoted to the study of singular backward stochastic Volterra integral equations. In Section 3.1, we discuss its well-posedness under quite weak singular conditions. In Section 3.2, we apply the aforementioned singular conditions to the forward stochastic Volterra system. In Section 3.3, we give several different examples to show the usefulness of the previous results.

\subsection{Well-posedness of singular BSVIEs}
In this section, we investigate the well-posedness of BSVIE \rf{Type-II}.
For readers' convenience, we rewrite it here,
\bel{BSVIE}
Y(t) = \psi(t) + \int_t^T g(t,s,Y(s),Z(t,s),Z(s,t))ds - \int_t^T Z(t,s)dW(s), \qq t\in \T,
\ee
where the free term $\psi(\cd) \in L^2_{\cF_T}(0,T;H)$ and the generator $g$ are given. A pair of $\lt(Y(\cd),Z(\cd,\cd)\rt) \in \cH^2[0,T]$ (Recalling the definition of $\cH^2[0,T]$ in Section 2) is referred to as an {\it adapted $M$-solution} of BSVIE \rf{BSVIE} if \rf{BSVIE} holds in the usual It\^o's sense for almost all $t\in \T$ and the following holds:
$$ Y(t) = \dbE Y(t) + \int_0^t Z(t,s) dW(s), \qq \ae \;\, t\in \T.$$

We introduce the following assumptions.\\

{\bf (H3.1)} Let $g: \D^* \times H \times \cL_2^0 \times \cL_2^0 \times \O \mapsto H$ be %$\cB(\D^* \times H \times \cL_2^0 \times \cL_2^0) \otimes \cF_T$-
measurable such that $s\mapsto g(t,s,y,z_1,z_2)$ is $\dbF$-progressively measurable for all $(t,y,z_1,z_2)\in \T \times H \times \cL_2^0 \times \cL_2^0$ and $g(t,s,0,0,0) =0$. Moreover, it holds that
$$\begin{aligned} |g(t,s,y,z_1,z_2) - g(t,s,\hat{y},\hat{z}_1,\hat{z}_2)|_H \les& L_y(t,s)|y-\hat{y}|_H + L_{z_1}(t,s)|z_1-\hat{z}_1|_{\cL_2^0} + L_{z_2}(t,s)|z_2-\hat{z}_2|_{\cL_2^0}, \\
	&\qq\qq\qq \forall(t,s,y,z_1,z_2),(t,s,\hat{y},\hat{z}_1,\hat{z}_2) \in \D^* \times H \times \cL_2^0 \times \cL_2^0,
\end{aligned}$$
where $L_y \in L^2(\D^*), L_{z_2} \in \sL^2(\D^*)$ and $\displaystyle \sup_{t\in (0,T)} \int_t^T L_{z_1}(t,s)^2ds <\i$.\\

\rm

\br{rmk41} In both Anh et al \cite{Y 2011} and Yong \cite{Y 2008}, the authors made similar Lipschitz conditions as that in (H3.1) and assumed that all the coeffcients satisfy simultaneously that
%(
\bel{2+eps} \sup_{t\in [0,T]} \int_t^T f(t,s)^{2+\e}ds < \i,\q f\= L_y, L_{z_1}, L_{z_2}, \text{ for some } \e >0. \ee

In contrast, we give a deeper exploration of the singularities for the three coefficients, and impose appropriate assumptions among them. In fact, the assumption of $L_y$ is the weakest while that of $L_{z_2}$ is the strongest. The condition of $L_{z_1}$ is in between the previous two cases. We point out that even for the
strongest case, the following example shows that it is weaker than that in Anh et al \cite{Y 2011} and Yong \cite{Y 2008}.
Let
$$h(s) \= \frac{1}{\sqrt{s} \ln s},\q s\in (0,T).$$
After some basic calculations, we can verify $\ds h(\cd) \in L^2(0,T;\dbR_+)\setminus  \bigcup_{\e>0} L^{2+\e}(0,T;\dbR_+).$ Let
$$L_y(t,s)= L_{z_1}(t,s)= L_{z_2}(t,s) \= h(s-t), \q (t,s) \in \D^*.$$
Then $f(\cd,\cd) \in \sL^2(\D^*) \subset L^2(\D^*)$ for $f\= L_y, L_{z_1}, L_{z_2}$, but \rf{2+eps} does not hold anymore.
\er

\br{} The case of $g(t,s,0,0,0) \neq 0$ can be treated as follows. Define a new free term
$$ \ti{\psi}(t) \= \psi(t) + \int_t^T g(t,s,0,0,0)ds, \q t\in [0,T].$$
We can prove $\ti{\psi}(\cd) \in L^2_{\cF_T}(0,T;H)$ if $\displaystyle \dbE\int_0^T \(\int_t^T |g(t,s,0,0,0)|_H ds\)^2 dt <\i$ is provided.
\er

%\bde{adapted M}

%\ede

Now we give the main result in this section.

\bt{43} Let {\bf (H3.1)} hold. Then for any $\psi(\cd) \in L^2_{\cF_T}(0,T;H)$, \rf{BSVIE} admits a unique adapted M-solution in $\cH^2[0,T]$. Moreover, it holds that
$$
\dbE \lt\{\int_{0}^{T} |Y(t)|_H^2dt  + \int_0^{T} \int_{0}^{T} |Z(t,s)|_{\cL_2^0}^2 ds dt\rt\} \les C \dbE \int_{0}^{T} |\psi(t)|_H^2 dt.
$$
Let $\bar g$ also satisfy {\bf (H3.1)}, $\bar \psi(\cd) \in L^2_{\cF_T}(0,T;H)$,  and $(\bar Y(\cd), \bar Z(\cd,\cd)) \in \cH^2[0,T]$ be the adapted $M$-solution of \rf{BSVIE} with $g$ and $\psi(\cd)$ replaced by $\bar g$ and $\bar \psi(\cd)$, respectively. Then the following stability estimate holds,
$$\begin{aligned}
	&\dbE \lt\{\int_{0}^{T} |Y(t)-\bar Y(t)|_H^2dt  + \int_{0}^{T} \int_{0}^{T} |Z(t,s)-\bar Z(t,s)|_{\cL_2^0}^2 dsdt \rt\}\\
	%
	%\les & C\dbE \int_{0}^{T} |\hat\psi(t)|^2 dt \\
	%
	\les & C \Big\{ \dbE \int_{0}^{T} |\psi(t)- \bar \psi(t)|_H^2 dt \\
	&\qq + \dbE \int_{0}^{T} \( \int_t^T |g(t,s,\bar Y(s), \bar Z(t,s), \bar Z(s,t)) - \bar g (t,s,\bar Y(s), \bar Z(t,s), \bar Z(s,t))|_Hds \)^2 dt \Big\}.
\end{aligned}$$
\et

\vskip 5mm

Before giving the proof of Theorem \ref{43}, we give some lemmas which will be uesful later. %and can be found in \cite{Y 2008}.
For any $R, S \in [0, T )$, we consider the following $H$-valued stochastic integral equation:
\bel{l}
\l(t,r) = \psi(t) + \int_r^T h(t,s,\mu(t,s))ds - \int_r^T \mu(t,s)dW(s), \q t\in [S,T],\; r\in[R,T],\ee
where $h:[S,T] \times [R,T] \times \cL_2^0 \times \O \mapsto H$ is given. The unknown process is $(\l(\cd , \cd ),\mu(\cd  , \cd ))$, for which $(\l(t,\cd ),\mu(t,\cd ))$ is $\dbF$-adapted for all $t \in [S, T ]$. We
may regard the above as

\vskip 2mm
\no $\bullet$ a family of {infinite dimensional} BSDEs on $[R, T ]$, parameterized by $t \in [S, T ]$;

\vskip 2mm
\no $\bullet$  a family of {infinite dimensional} stochastic Fredholm{-type} integral equations (SFIEs, for short) on $[S, T ]$, parameterized by $r\in [R, T ]$.\\

We introduce the following assumptions concerning the generator $h$ of \rf{l}.\\

{\bf (H3.2)} Let $R,S \in [0,T)$, and $h:[S,T] \times [R,T] \times \cL_2^0 \times \O \mapsto H$ be %$\cB([S,T] \times [R,T] \times \cL_2^0) \otimes \cF_T$-
measurable such that $s \mapsto h(t,s,z)$ is $\dbF$-progressively measurable for all $(t,z) \in [S,T] \times \cL^0_2 $ and
$$ \int_S^T \dbE \(\int_R^T |h(t,s,0)|_H ds\)^2dt < \i.$$
Moreover, the following holds:
$$ |h(t,s,z)- h(t,s,\bar{z})|_H \les L(t,s) |z-\bar{z}|_{\cL_2^0},\q (t,s) \in [S,T] \times [R,T], \; z,\bar{z} \in \cL^0_2, \; \as,$$
where $L:[S,T] \times [R,T] \to [0,\i)$ is a determinitic function such that
$$ \sup_{t\in [S,T]} \int_R^T L(t,s)^2 ds < \i.$$

\vskip 2mm
\bl{4.1.}
Let {\bf (H3.2)} hold. Then for any $\psi(\cd) \in L_{\cF_T}^2(S,T;H)$,  \rf{l} %regarded as an infinite dimensional BSDE on $[R,T]$,
admits a unique adapted solution
$$(\l(t,\cd),\mu(t,\cd)) \in L_\dbF^2(\O;C([R,T];H)) \times L_\dbF^2(R,T;\cL_2^0)$$
for almost all $t\in [S,T]$. Moreover, if $\bar{h}$ also satisfies {\bf (H3.2)}, $\bar{\psi}(\cd) \in L_{\cF_T}^2(S,T;H)$, and $$(\bar{\l}(t,\cd), \bar{\mu}(t,\cd)) \in L_\dbF^2(\O;C([R,T];H)) \times L_\dbF^2(R,T;\cL_2^0)$$ is the unique adapted solution of \rf{l} with $(h,\psi)$ replaced by $(\bar{h},\bar{\psi})$. Then
$$ \begin{aligned}
	&\dbE \lt\{\sup_{r\in [R,T]}|\l(t,r) - \bar{\l}(t,r)|_H^2 + \int_R^T |\mu(t,s) - \bar\mu(t,s)|_{\cL_2^0}^2 ds\rt\} \\
	\les\;& C \dbE \lt\{ |\psi(t) - \bar{\psi}(t)|_H^2 + \(\int_R^T |h(t,s,\mu(t,s)) - \bar{h}(t,s,\mu(t,s))|_H ds\)^2
	\rt\}, \q \ae \; t\in [S,T].
\end{aligned}$$
\el
{\pf Since the proof is similar to the finite dimensional case in \cite{Y 2008}, we omit here. \ef}

\vskip 2mm
Next, let us look at two special cases of the above result. First, let $r=S \in [R,T)$ be fixed. Define
$$ \psi^S(t) \= \l(t,S), \q Z(t,s) \=\mu(t,s),\q t\in[R,S],\;s\in [S,T].$$
Then equation \rf{l} reads:
\bel{SFIE}
\psi^S(t) = \psi(t) + \int_S^T h(t,s,Z(t,s))ds - \int_S^T Z(t,s)dW(s), \q t\in[R,S].\ee
This is a stochastic Fredholm type integral equation valued in Hilbert spaces. A pair $$(\psi^S(\cd),Z(\cd,\cd)) \in L_{\cF_S}^2(R,S;H) \times L^2(R,S;L_\dbF^2(S,T;\cL_2^0))$$ satisfying
\rf{SFIE} in the usual It\^o sense is called an adapted solution of \rf{SFIE}. We note
that $\psi^S(t)$ is (only) required to be $\cF_S$-measurable for almost all $t \in [R, S]$, instead of $\dbF$-adaptiveness. According to Lemma \ref{4.1.}, we have the following result.
\bc{cor1}
Let {\bf (H3.2)} hold. Then for any $\psi(\cd) \in L_{\cF_T}^2(R,S;H)$, SFIE \rf{SFIE} admits a unique adapted solution $(\psi^S(\cd),Z(\cd,\cd)) \in L_{\cF_S}^2(R,S;H) \times L^2(R,S;L_\dbF^2(S,T;\cL_2^0))$.
\ec

\vskip 2mm
The second special case of \rf{l} is the following: Let $R=S$, and define
\be \notag \lt\{\begin{aligned}
	& Y(t) \= \l(t,t), \qq t\in[S,T],\\
    &Z(t,s) \= \mu(t,s), \q (t,s) \in \D^*[S,T].
\end{aligned} \rt.\ee
Then \rf{l} reads:
\bel{s}
Y(t) = \psi(t) + \int_t^T h(t,s,Z(t,s))ds - \int_t^T Z(t,s)dW(s), \q t\in [S,T].
\ee
This can be seen as a special case of  BSVIE \rf{BSVIE} in which the generator $h$ is independent of $Y(s)$ and $Z(s,t)$. We can define $Z(t, s)$ for $(t, s) \in \D[S, T ]$ by the martingale representation theorem:
$$Y(t) = \dbE[Y(t) |\cF_S] + \int_S^t Z(t,s)dW(s), \q t\in [S,T].$$
We have the following result.
\bc{cor2}
Let {\bf (H3.2)} hold. Then for any $\psi(\cd) \in L_{\cF_T}^2(S,T;H)$, BSVIE \rf{s} admits a unique adapted $M$-solution
$$(Y(\cd),Z(\cd,\cd)) \in L_{\dbF}^2(S,T;H) \times L^2(S,T;L_\dbF^2(S,T;\cL_2^0)).$$ Moreover, if $\bar{h}$ also satisfies {\bf (H3.2)}, $\bar{\psi}(\cd) \in L_{\cF_T}^2(S,T;H)$, and $(\bar{Y}(\cd), \bar{Z}(\cd,\cd)) \in L_{\dbF}^2(S,T;H) \times L^2(S,T;L_\dbF^2(S,T;\cL_2^0))$ is the unique adapted $M$-solution of \rf{s} with $(h,\psi)$ replaced by $(\bar{h},\bar{\psi})$, then
$$ \begin{aligned}
	&\dbE \lt\{|Y(t) - \bar{Y}(t)|_H^2 + \int_S^T |Z(t,s) - \bar{Z}(t,s)|_{\cL_2^0}^2 ds\rt\} \\
	\les\;& C \dbE \lt\{ |\psi(t) - \bar{\psi}(t)|_H^2 + \(\int_t^T |h(t,s,Z(t,s)) - \bar{h}(t,s,Z(t,s))|_H ds\)^2
	\rt\}, \q t\in [S,T].
\end{aligned}$$
\ec

\color{black}

With the above preparations, we are ready to give the proof of Theorem \ref{43}.

\begin{proof}[Proof of Theorem \ref{43}]
	
First we define $\cM^2\T $ as the space of all $(y(\cd),z(\cd,\cd)) \in \cH^2\T$ such that
$$ y(t) = \dbE y(t) + \int_0^t z(t,s) dW(s), \qq t\in \T.$$
We introduce an equivalent norm for $\cM^2\T$ as follows:
$$\lt\| (y(\cd),z(\cd,\cd)) \rt\|_{\cM^2[0,T]} \= \[ \dbE \int_0^T |y(t)|_H^2 dt + \dbE \int_0^T \(\int_t^T |z(t,s)|_{\cL_2^0}^2 ds\)dt\]^\frac{1}{2}.$$

Step 1: Fix a $(y(\cd),z(\cd,\cd)) \in \cM^2[S,T]$ with $S$ undertermined, let us consider the following equation:
\bel{51} Y(t) = \psi(t) + \int_t^T g(t,s,y(s),Z(t,s),z(s,t))ds - \int_t^T Z(t,s)dW(s), \qq t\in [S,T], \ee
for any $\psi(\cd) \in L^2_{\cF_T}(S,T;H)$. By Corollary \ref{cor2}, we can see that \rf{51} admits a unique adapted $M$-solution  $(Y(\cd),Z(\cd,\cd))$, and
\bel{0es}\begin{aligned}
	&\dbE \lt\{ |Y(t)|_H^2 + \int_t^T |Z(t,s)|_{\cL_2^0}^2 ds \rt\} \\
	\les & C \dbE \lt\{ |\psi(t)|_H^2 + \(\int_t^T |g(t,s,y(s),0,z(s,t))|_Hds\)^2 \rt\} \\
	\les & C \dbE \lt\{ |\psi(t)|_H^2 + \(\int_t^T \[L_y(t,s)|y(s)|_H + L_{z_2}(t,s) |z(s,t)|_{\cL_2^0}\] ds\)^2 \rt\} \\
	\les &C \dbE \lt\{ |\psi(t)|_H^2 + \int_t^T L_y(t,s)^2 ds \cd \int_t^T |y(s)|_H^2 ds + \int_t^T L_{z_2}(t,s)^2 ds \cd \int_t^T |z(s,t)|_{\cL_2^0}^2 ds \rt\}.
\end{aligned}\ee
By integrating on $[S,T]$ w.r.t. $t$, we derive that
\bel{1es} \begin{aligned}
	& \dbE \lt\{\int_S^T |Y(t)|_H^2 dt + \int_S^T \int_t^T |Z(t,s)|_{\cL_2^0}^2 dsdt \rt\} \equiv \lt\| (Y(\cd),Z(\cd,\cd)) \rt\|_{\cM^2[S,T]}^2 \\
	\les C & \dbE \[\int_S^T |\psi(t)|_H^2 dt + \int_S^T \int_t^T L_y(t,s)^2 dsdt \cd \int_S^T |y(t)|_H^2 dt + \sup_{t\in [S,T]} \int_t^T L_{z_2}(t,s)^2 ds \cd \int_S^T\int_t^T |z(s,t)|_{\cL_2^0}^2 dsdt
	\] \\
	\les C & \( \|L_y(\cd,\cd)\|^2_{L^2(\D^*[S,T])} + \|L_{z_2}(\cd,\cd)\|^2_{\sL^2(\D^*[S,T])}
	\) \cd \[\dbE \int_S^T |\psi(t)|_H^2 dt + \| (y(\cd),z(\cd,\cd)) \|_{\cM^2[S,T]}^2\].
\end{aligned} \ee

\no Thus we define a map $\Th: \cM^2[S,T] \mapsto \cM^2[S,T]$ by
$$ \Th(y(\cd),z(\cd,\cd)) = (Y(\cd),Z(\cd,\cd)), \qq\forall (y(\cd),z(\cd,\cd)) \in \cM^2[S,T].$$
Now we show that the map $\Th$ is contractive for some $0 \les S <T$. Indeed, for another $(\bar{y}(\cd),\bar{z}(\cd,\cd)) \in \cM^2[S,T]$ and $\Th(\bar{y}(\cd),\bar{z}(\cd,\cd)) = (\bar{Y}(\cd),\bar{Z}(\cd,\cd))$. By the stability estimate in Corollary \ref{cor2},
\bel{1contractive}\begin{aligned}
  &\dbE \int_S^T |Y(t)-\bar{Y}(t)|_H^2 dt + \dbE \int_S^T \int_t^T |Z(t,s) - \bar{Z}(t,s)|_{\cL_2^0}^2ds dt \\
  \les\; & C\dbE \int_S^T \[ \int_t^T |g(t,s,y(s),Z(t,s),z(s,t)) - g(t,s,\bar{y}(s),Z(t,s),\bar{z}(s,t))|_Hds \]^2 dt \\
  \les\; & C\dbE \int_S^T \[ \int_t^T L_y(t,s)|y(s) - \bar{y}(s)|_Hds\]^2dt + C\dbE \int_S^T \[ \int_t^T L_{z_2}(t,s)|z(s,t) - \bar{z}(s,t)|_{\cL_2^0}ds\]^2dt \\
  \les\; & C\dbE \int_S^T \[ \int_t^T L_y(t,s)^2 ds\cd \int_t^T |y(s) - \bar{y}(s)|_H^2 ds \]dt \\
  &\qq\qq\qq \qq +  C\dbE \int_S^T \[ \int_t^T L_{z_2}(t,s)^2 ds\cd \int_t^T |z(s,t) - \bar{z}(s,t)|_{\cL_2^0}^2 ds \]dt \\
  \les\; & C\int_S^T\int_t^T L_y(t,s)^2 dsdt \cd \dbE \int_S^T |y(s) - \bar{y}(s)|_H^2 ds \\
  &\qq\qq\qq\qq + C\sup_{t\in (S,T)}\int_t^T L_{z_2}(t,s)^2 ds \cd \dbE \int_S^T \int_S^t |z(t,s) - \bar{z}(t,s)|_{\cL_2^0}^2 dsdt \\
  \les\; & C\[ \int_S^T\int_t^T L_y(t,s)^2 dsdt + \sup_{t\in (S,T)}\int_t^T L_{z_2}(t,s)^2 ds \] \cd \dbE \int_S^T |y(t) - \bar{y}(t)|_H^2 ds.
\end{aligned}\ee
Recall that $L_y \in L^2(\D^*)$ and $L_{z_2} \in \sL^2(\D^*)$. {With the above $C$,} there exists a partition $\{U_i\}_{i=0}^m$ of $[0,T]$ with $0=U_0 < U_1 <\cds <U_m =T$ such that
\bel{bpart} C \[\int_{U_i}^{U_{i+1}} \int_{t}^{U_{i+1}} L_y(t,s)^2 dsdt + \sup_{ t\in (U_i,U_{i+1})} \int_t^{U_{i+1}} L_{z_2}(t,s)^2ds \] \les \frac{1}{2}, \q \forall i=0,1,\cds,m-1.\ee
In other words, if we choose $S= U_{m-1}$, then
$$ C\[ \int_S^T\int_t^T L_y(t,s)^2 dsdt + \sup_{t\in (S,T)}\int_t^T L_{z_2}(t,s)^2 ds\] \les \frac{1}{2}. $$
Therefore, the map $\Th$ admits a unique fixed point $(Y(\cd),Z(\cd,\cd))\in \cM^2[S,T]$ which is the unique solution of \rf{BSVIE} over $[U_{m-1},T]$. This step determines the values $(Y(t),Z(t,s))$ for $(t,s) \in [U_{m-1},T]\times [U_{m-1},T]$. By \rf{1es} and \rf{bpart}, we have the following estimate
\bel{1es2}
\dbE \lt\{\int_{U_{m-1}}^T |Y(t)|_H^2 dt + \int_{U_{m-1}}^T \int_t^T |Z(t,s)|_{\cL_2^0}^2 dsdt \rt\} \les C \dbE \int_{U_{m-1}}^T |\psi(t)|_H^2 dt.
\ee

Step 2: We determine the values $Z(t,s)$ of $Z(\cd,\cd)$ for $(t,s)\in [U_{m-1}, T]\times[U_{m-2}, U_{m-1}]$ by means of martingale representation theorem, i.e.
$$\dbE[Y(t)|\cF_{U_{m-1}}] = \dbE[Y(t)|\cF_{U_{m-2}}] + \int_{U_{m-2}}^{U_{m-1}} Z(t,s)dW(s), \qq t\in
[U_{m-1}, T].$$
So
$$ \dbE \int_{U_{m-2}}^{U_{m-1}} |Z(t,s)|_{\cL_2^0}^2ds \les \dbE |Y(t)|_H^2, \q t\in [U_{m-1},T].$$
By \rf{1es2} and integrating the above inequality, we have
\bel{2es}
\dbE \int_{U_{m-1}}^T \int_{U_{m-2}}^{U_{m-1}} |Z(t,s)|_{\cL_2^0}^2dsdt \les \dbE \int_{U_{m-1}}^T |Y(t)|_H^2dt \les C\dbE \int_{U_{m-1}}^T |\psi(t)|_H^2 dt.
\ee

Recall that, we have determined $(Y(t),Z(t,s))$ for $(t,s) \in [U_{m-1},T]\times [U_{m-2},T]$ up to now. Combining \rf{1es2} and \rf{2es}, we get
\bel{2es2}
\dbE \lt\{\int_{U_{m-1}}^T |Y(t)|_H^2 dt + \int_{U_{m-1}}^T \int_{U_{m-2}}^T |Z(t,s)|_{\cL_2^0}^2 dsdt \rt\} \les C \dbE \int_{U_{m-1}}^T |\psi(t)|_H^2 dt.
\ee

Step 3: For $(t,s) \in
[U_{m-2}, U_{m-1}] \times
[U_{m-1}, T]$, we know that the values $Y(s)$ and $Z(s,t)$ are determined by Step 2. In this step, we determine the value of $Z(t,s)$ for $(t,s) \in
[U_{m-2}, U_{m-1}] \times
[U_{m-1}, T]$ by solving a stochastic Fredholm integral equation, i.e.
\bel{52} \psi^{U_{m-1}}(t) = \psi(t) + \int_{U_{m-1}}^T g^{U_{m-1}}(t,s,Z(t,s))ds - \int_{U_{m-1}}^T Z(t,s)dW(s), \qq t\in[U_{m-2},{U_{m-1}}],\ee
where $$g^{U_{m-1}}(t,s,z)\= g(t,s,Y(s),z,Z(s,t)),\q (t,s,z)\in [U_{m-2}, U_{m-1}] \times [U_{m-1} ,T]\times \cL_2^0.$$
Thanks to Corollary \ref{cor1}, \rf{52} has a unique adapted solution $$(\psi^{U_{m-1}}(\cd),Z(\cd,\cd)) \in L_{\cF_{U_{m-1}}}^2(U_{m-2},{U_{m-1}};H) \times L^2(U_{m-2},{U_{m-1}};L_\dbF^2({U_{m-1}},T;\cL_2^0)),$$ and similar to \rf{1es}, it holds that
$$\begin{aligned}
	&\dbE \lt\{ |\psi^{U_{m-1}}|_H^2 + \int_{U_{m-1}}^T |Z(t,s)|_{\cL_2^0}^2 ds \rt\} \\
	\les & C \dbE \lt\{ |\psi(t)|_H^2 + \(\int_{U_{m-1}}^T |g(t,s,Y(s),0,Z(s,t))|_Hds\)^2 \rt\} \\
	\les &C \dbE \lt\{ |\psi(t)|^2 + \int_{U_{m-1}}^T L_y(t,s)^2 ds \cd \int_{U_{m-1}}^T |Y(s)|_H^2 ds + \int_{U_{m-1}}^T L_{z_2}(t,s)^2 ds \cd \int_{U_{m-1}}^T |Z(s,t)|_{\cL_2^0}^2 ds \rt\}.
\end{aligned}$$
Therefore,
\bel{3es}\begin{aligned}
	&\dbE \int_{U_{m-2}}^{U_{m-1}} \lt\{ |\psi^{U_{m-1}}|_H^2 + \int_{U_{m-1}}^T |Z(t,s)|_{\cL_2^0}^2 ds   \rt\}dt\\
    \les & C \dbE \lt\{  \int_{U_{m-2}}^{U_{m-1}} |\psi(t)|_H^2 dt + \int_{U_{m-1}}^T |Y(t)|_H^2 dt + \int_{U_{m-2}}^{U_{m-1}} \int_{U_{m-1}}^T |Z(s,t)|_{\cL_2^0}^2 dsdt \rt\}
    \les  C \dbE \int_{U_{m-2}}^{T} |\psi(t)|_H^2 dt.
\end{aligned}\ee
The last inequality holds due to \rf{2es2}.\\

Step 4: So far, we have uniquely determined
\be \left\{ \begin{aligned}
	&Y(t), \qq t\in[U_{m-1}, T], \\
	&Z(t,s), \qq (t,s) \in ([U_{m-1}, T]\times [U_{m-2},T]) \cup ([U_{m-2},U_{m-1}]\times [U_{m-1},T]).
\end{aligned}\right.\ee
Now, we consider
\bel{53} Y(t) = \psi^{U_{m-1}}(t) + \int_t^{U_{m-1}} g(t,s,Y(s),Z(t,s),Z(s,t))ds - \int_t^{U_{m-1}} Z(t,s)dW(s), \qq t \in [U_{m-2}, U_{m-1}].\ee
Since $\psi^{U_{m-1}}(\cd)$ is  ${\cF_{U_{m-1}}}$-measurable, \rf{53} is a BSVIE over $[{U_{m-2}},{U_{m-1}}]$. Hence, recalling the inequality \rf {bpart} (by choosing $U_i= U_{m-2}$ and $U_{i+1}=U_{m-1}$) and by the similar procedures in Step 1, we are able to show that \rf{53} is solvable on $[{U_{m-2}},{U_{m-1}}]$. Moreover,  by \rf{3es}, it holds that
\bel{4es}\begin{aligned}
	&\dbE \lt\{\int_{U_{m-2}}^{U_{m-1}} |Y(t)|_H^2dt  + \int_{U_{m-2}}^{U_{m-1}} \int_{U_{m-2}}^{U_{m-1}} |Z(t,s)|_{\cL_2^0}^2 dsdt \rt\} \\
    \les & C\dbE \int_{U_{m-2}}^{U_{m-1}} |\psi^{U_{m-1}}(t)|_H^2 dt
    \les   C \dbE \int_{U_{m-2}}^{T} |\psi(t)|_H^2 dt.
\end{aligned}\ee

This solvability determines $(Y(t),Z(t,s))$ for $(t,s)\in [{U_{m-2}},{U_{m-1}}]\times [{U_{m-2}},{U_{m-1}}]$. As a consequence, we obtain the unique solvability of BSVIE \rf{BSVIE} on $[U_{m-2},T]$ and  the following estimate holds by combining \rf{2es2}, \rf{3es}, \rf{4es},
$$
\dbE \lt\{\int_{U_{m-2}}^{T} |Y(t)|_H^2dt  + \int_{U_{m-2}}^{T} \int_{U_{m-2}}^{T} |Z(t,s)|_{\cL_2^0}^2 ds dt\rt\} \les C \dbE \int_{U_{m-2}}^{T} |\psi(t)|_H^2 dt.
$$

By the similar procedures, we can complete the proof by induction.\\

At last, we give the stability estimate. Suppose that $(Y(\cd),Z(\cd,\cd))$ and $(\bar Y(\cd),\bar Z(\cd,\cd))$ are adapted $M$-solution of \rf{BSVIE} corresponding to $(g,\psi)$ and $(\bar g,\bar\psi)$ respectively. Denote
$$\hat Y(t) \= Y(t)- \bar Y(t), \q \hat Z(t,s) \= Z(t,s)- \bar Z(t,s).
$$
\no Observe that $\hat Y(\cd)$ satisfies the following BSVIE,
$$
\hat Y(t) = \hat \psi(t) + \int_t^T \hat g(t,s,\hat Y(s), \hat Z(t,s), \hat Z(s,t)) ds -\int_t^T \hat Z(t,s) dW(s),
$$
where
$$\left\{\begin{aligned}
&\hat \psi(t) \= \psi(t)- \bar \psi(t) + \int_t^T [g(t,s,\bar Y(s), \bar Z(t,s), \bar Z(s,t)) - \bar g (t,s,\bar Y(s), \bar Z(t,s), \bar Z(s,t))] ds,\\
&\hat g(t,s,y,z_1,z_2) \= g(t,s, y+ \bar Y(s), z_1+ \bar Z(t,s),z_2+ \bar Z(s,t)) - g(t,s,\bar Y(s), \bar Z(t,s), \bar Z(s,t)).\end{aligned}\rt.$$
It's easy to check that the generator $\hat g$ satisfies
assumption {\bf (H3.1)} and the following stability estimate holds,
$$\begin{aligned}
	&\dbE \lt\{\int_{0}^{T} |\hat Y(t)|_H^2dt  + \int_{0}^{T} \int_{0}^{T} |\hat Z(t,s)|_{\cL_2^0}^2 ds dt\rt\}
	\les  C\dbE \int_{0}^{T} |\hat\psi(t)|_H^2 dt \\
	\les & C \Big\{ \dbE \int_{0}^{T} |\psi(t)- \bar \psi(t)|_H^2 dt \\
	&\qq + \dbE \int_{0}^{T} \( \int_t^T |g(t,s,\bar Y(s), \bar Z(t,s), \bar Z(s,t)) - \bar g (t,s,\bar Y(s), \bar Z(t,s), \bar Z(s,t))|_Hds \)^2 dt \Big\}.
\end{aligned}$$
\ef

\br{remark37} As to the above proof, we have the following remarks to point out.
\begin{itemize}
	\item[(i)]
 In contrast with \cite{Y 2008}, it is worthy mentioning some technology differences due to the coefficients' stronger singularity in the current paper (see e.g. Remark \ref{rmk41}). While the author in \cite{Y 2008} estimated the Lips. terms in the same way by H\"older inequality, here we considerably relax the coefficients $L_y, L_{z_1}, L_{z_2}$ and treat them differently. In fact, for the $L_{z_1}$, we rely on the BSDEs result, while for $L_y$, $L_{z_2}$, we use the singularities similar to that in Hamaguchi--Wang \cite{Hamaguchi}.

 \item[(ii)] The proof of Theorem \ref{43} seems to be lengthy and technical.
 One may consider it from a new point by defining an equivalent norm on the solution space as follows (see also \cite{W 2012} for instance),
 $$\lt\| (y(\cd),z(\cd,\cd)) \rt\|_{\cM^{2,\b}[0,T]} \= \[ \dbE \int_0^T e^{2\b t} |y(t)|_H^2 dt + \dbE \int_0^T e^{\b t} \(\int_t^T |z(t,s)|_{\cL_2^0}^2 ds\)dt\]^\frac{1}{2}.$$
 However, in this case, the assumptions on $L_y, L_{z_2}$  need to be strengthened, i.e., there exists some $\e > 0$, such that
 $$ \int_0^T \( \int_t^T L_y(t,s)^{2(1+\e)} ds \)^\frac{1}{1 +\e} dt + \sup_{ t\in [0,T]} \int_t^T L_{z_2}(t,s)^{2(1+\e)} ds < \i.
 $$
 This would end up with the essentially same singularity as that in \cite{Y 2008}.

\item[(iii)]
In \cite{Hamaguchi1}, the author considered an infinite horizon backward stochastic Volterra integral equations with convolution-type singular kernel.
Even when their system reduces to the finite interval case, it seems that their developed method is not well applicable to our scenario.
\end{itemize} \er

%{\color{red} Is it possible to incorporate the above remarks together?}

\subsection{Well-posedness of singular forward SVIEs}

In Section 3.1, we have proved the well-posedness of the singular backward stochastic Volterra system \rf{BSVIE} in the sense of adapted M-solution. It is found that the aforementioned singular assumption can be well fit into the forward system as follows:
\bel{fSVIE}
X(t) = \f(t) + \inte A(t,s,X(s))ds + \inte B(t,s,X(s))dW(s), \q t\in \T.
\ee
In this part, we present the solvability of the above \rf{fSVIE} under the following  assumptions.\\

{\bf (H3.3)} Let maps $A: \D \times H \times \Omega \mapsto H$, $B: \D \times H \times \Omega \mapsto \cL_2^0$ be measurable. For each $(t,x) \in \T \times H, s \mapsto \lt(A(t,s,x), B(t,s,x)\rt)$ is $\dbF$-adapted on $[0,t]$. And,
$$A(t,s,0) =0, \text{ and }B(t,s,0)=0 \text{\q  for \ae\; $(t,s) \in \D, \;\as$} $$
Moreover,
there are $K_1 \in L^2(\D)$ and $K_2 \in \mathscr{L}^2(\D)$ such that for any $(t,s) \in \D%= \{(t,s) \in [0,T]^2 \;|\; 0\les s < t \les T\}
$ and $x, y \in H$, it holds that
\be \notag \begin{aligned}
	|A(t,s,x) - A(t,s,y)|_H &\les K_1(t,s) |x-y|_H, \\
	|B(t,s,x) - B(t,s,y)|_{\cL_2^0} &\les K_2(t,s) |x-y|_H.
\end{aligned} \ee

\br{nonzero} The case with nonzero $A(t,s,0)$ and $B(t,s,0)$ can be treated as follows. As a matter of fact, suppose that $A$ and $B$ satisfy the Lipschitz condition in Assumption {\bf (H3.1)}, and
$$ \dbE \int_0^T \(\int_0^t |A(t,s,0)|_H ds\)^2dt < \i, \q \dbE \int_0^T \int_0^t |B(t,s,0)|_{\cL_2^0}^2dsdt<\i. $$
%
%\tb{ $\dbE$ may be not needed! }
%
We can rewrite \rf{fSVIE} as follows:
\be \notag \begin{aligned}
	X(t) =&\f(t) + \int_0^t A(t,s,0)ds + \int_0^t B(t,s,0)dW(s) \\
	&\qq \qq+ \int_0^t \tilde{A}(t,s,X(s))ds + \int_0^t \tilde{B}(t,s,X(s))dW(s), \q t\in \T,
\end{aligned} \ee
where
$$\tilde{A}(t,s,x) \= A(t,s,x) - A(t,s,0), \q \tilde{B}(t,s,x) \= B(t,s,x) - B(t,s,0).$$
Clearly,
$$ \tilde{A}(t,s,0)=0, \q \tilde{B}(t,s,0)=0,$$
and $\wt{A}, \wt{B}$ satisfy the assumption {\bf (H3.3)}.
Let
$$\tilde{\f}(t) \= \f(t) + \int_0^t A(t,s,0)ds + \int_0^t B(t,s,0)dW(s), \q t\in \T,$$
be the new free term. Moreover, we can prove $\tilde{\f}(\cd) \in L_\dbF^2(0,T;H)$ if $\f(\cd) \in L_\dbF^2(0,T;H)$.\\
\er

\bt{thm3} Let {\bf (H3.3)} hold. Then for any $\f(\cd) \in L_\dbF^2(0,T;H)$, the equation \rf{fSVIE} admits a unique solution $X(\cd) \in L_\dbF^2(0,T;H)$. Moreover, the following estimate holds:
\bel{estimate} \lt\| X(\cd) \rt\|_{L_\dbF^2(0,T;H)} \les  C \lt\| \f(\cd) \rt\|_{L_\dbF^2(0,T;H)}.\ee

\no Let $A'$ and $B'$ satisfy {\bf (H3.3)},  $\f'(\cd) \in L_\dbF^2(0,T;H)$, and  $X'(\cd) \in L_\dbF^2(0,T;H)$ be the solution of SVIE \eqref{fSVIE} corresponding to $(\f', A', B')$. Then it holds that
\begin{equation}\notag
	\begin{aligned}
		& \[\dbE \int_0^T |X(t) - X'(t)|_H^2 dt \]^\frac{1}{2} \\
		\les & C \Big\{ \dbE \int_0^T \[ |\f(t) - \f'(t)|_H^2 + \(\int_0^t |A(t,s,X'(s)) - A'(t,s,X'(s))|_H ds\)^2 \\
		& \qq\qq\qq\qq \qq\qq +  \int_0^t |B(t,s,X'(s)) - B'(t,s,X'(s))|_{\cL_2^0}^2 ds \]dt\Big\}^\frac{1}{2}.
	\end{aligned}
\end{equation}

\et

\pf We use the Banach fixed point theorem to give the existence and uniqueness of the adapted solution of \rf{fSVIE}.

Step 1:  For any $S\in (0,T]$ and
$x(\cd) \in L_\dbF^2(0,S;H)$, we define
$$
\cS[x(\cd)](t) \= \f(t) + \int_0^t A(t,s,x(s))ds + \int_0^t B(t,s,x(s))dW(s), \q t\in [0,S].
$$
%F
Thus,
\bel{mapS} \begin{aligned}
	&\dbE \int_0^S |\cS[x(\cd)](t)|_H^2 dt \\
	\les C&\Big[ \dbE \int_0^S |\f(t)|_H^2 dt + \dbE \int_0^S \( \int_0^t K_1(t,s)|x(s)|_Hds\)^2 dt
	+\dbE \int_0^S \int_0^t K_2(t,s)^2|x(s)|_H^2dsdt \Big] \\
	\les  C&\Big[ \dbE \int_0^S |\f(t)|_H^2 dt + \int_0^S \int_0^t K_1(t,s)^2dsdt \cd \dbE \int_0^S |x(s)|_H^2ds + \dbE \int_0^S \int_s^S K_2(t,s)^2dt \cd |x(s)|_H^2ds \Big]\\
	\les C&\Big[ \dbE \int_0^S |\f(t)|_H^2 dt  + \lt(\|K_1\|_{L^2(\D)} + \|K_2\|_{\sL^2(\D)}\rt)  \cd  \dbE\int_0^S |x(t)|_H^2 dt \Big].
\end{aligned}\ee
Hence, $\cS$ maps $L_\dbF^2(0,S;H)$ to itself.\\

Step 2: Given $T_1 \in (0,T]$ to be determined later.
We prove that the above defined $\cS$ is contractive on the interval $[0,T_1]$.
For any $x(\cd), y(\cd) \in L_\dbF^2(0,T_1;H)$, %(with $T_1$ undetermined),similar to the above calculations
thanks to the H\"older inequality %BDG inequality%
and Fubini theorem, we have
$$ \begin{aligned}
	& \lt\| \cS[x(\cd)](t) - \cS[y(\cd)](t) \rt\|_{L_\dbF^2(0,T_1;H)}^2 \\
	\les & 2\dbE \int_0^{T_1} \( \int_0^t K_1(t,s)|x(s)-y(s)|_H ds\)^2 dt
	+2\dbE \int_0^{T_1} \int_0^t K_2(t,s)^2|x(s)-y(s)|_H^2dsdt \\
	\les & 2\int_0^{T_1} \int_0^t K_1(t,s)^2dsdt \cd \dbE \int_0^{T_1} |x(s)-y(s)|_H^2ds + 2\dbE \int_0^{T_1} \int_s^{T_1} K_2(t,s)^2dt \cd |x(s)-y(s)|_H^2ds. \\
\end{aligned}$$
By the definitions of $K_1, K_2$, there exists a partition $\{T_i\}_{i=0}^m$ of $[0,T]$ with $0=T_0 < T_1 <\cds <T_m =T$ such that
\bel{K1K2} \int_{T_i}^{T_{i+1}} \int_{T_i}^{t} K_1(t,s)^2 dsdt + \sup_{ s\in (T_i,T_{i+1})} \int_s^{T_{i+1}} K_2(t,s)^2dt \les \frac{1}{4}, \q \forall i=0,1,\cds,m-1.\ee

\no In particular, take $i=0$, i.e.
$$\int_0^{T_1} \int_0^t K_1(t,s)^2dsdt + \sup_{ s\in (0,T_{1})} \int_s^{T_{1}} K_2(t,s)^2dt \les \frac{1}{4}.$$
Hence, $\cS$ is contractive on $L_\dbF^2(0,T_1;H)$. Therefore \rf{fSVIE} admits a unique solution  $X_1(\cd) \in L_\dbF^2(0,T_1;H)$ on $[0,T_1]$. Further, by \rf{mapS} and noting \rf{K1K2}, $X(\cd)$ satisfies the estimate \rf{estimate} with $T$ replaced by $T_1$.\\

Step 3:
Next, we rewrite \rf{fSVIE} as follows:
$$ X(t) = \hat\f(t) + \int_{T_1}^t A(t,s,X(s)) ds + \int_{T_1}^t B(t,s,X(s)) dW(s), \q t\in [T_1,T],$$
where
$$\hat\f(t) \= \f(t) + \int_0^{T_1} A(t,s,X_1(s))ds + \int_0^{T_1} B(t,s,X_1(s))dW(s), \q t\in [T_1,T].$$
It is straightforward to verify $\hat\f(\cd) \in L_\dbF^2(T_1,T;H)$. For any $x(\cd) \in  L_\dbF^2(T_1,T_2;H)$, define
$$
\cS[x(\cd)](t) \= \hat\f(t) + \int_{T_1}^t A(t,s,x(s))ds + \int_{T_1}^t B(t,s,x(s))dW(s), \q t\in [T_1,T_2].
$$
Then by similar procedures in Step 1--2 and noting \rf{K1K2}, we can obtain the unique solution of \rf{fSVIE} on the interval $[T_1,T_2]$.
By repeating the above procedures, we have the solution on the whole interval $[0,T]$. Moreover, the estimate \rf{estimate} holds.\\

Observe that $\bar{X}(\cd) \= X(\cd) - X'(\cd) \in L_\dbF^{2}(0,T;H) $ is the solution of the following SVIE
$$ \bar{X}(t)=\bar{\varphi}(t)+\int_{0}^{t} \bar{A}(t, s, \bar{X}(s)) \mathrm{d} s+\int_{0}^{t} \bar{B}(t, s, \bar{X}(s)) \mathrm{d} W(s),\q  t \in [0,T],$$
where
$$\begin{aligned}
	&\bar{\varphi}(t)\=\varphi(t)-\varphi^{\prime}(t)+\int_{0}^{t}\left\{A\left(t, s, X^{\prime}(s)\right)-A^{\prime}\left(t, s, X^{\prime}(s)\right)\right\} \mathrm{d} s+\int_{0}^{t}\left\{B\left(t, s, X^{\prime}(s)\right)-B^{\prime}\left(t, s, X^{\prime}(s)\right)\right\} \mathrm{d} W(s), \\
	& \bar{A}(t, s, x)\=A\left(t, s, x+X^{\prime}(s)\right)-A\left(t, s, X^{\prime}(s)\right), \q \bar{B}(t, s, x) \=B\left(t, s, x+X^{\prime}(s)\right)-B\left(t, s, X^{\prime}(s)\right).
\end{aligned}$$
Note that $\bar{\f} \in {L_\dbF^{2}(0,T;H)}$, $\bar{A}$ and $\bar{B}$ satisfy {\bf (H3.3)} with the Lipschitz functions $K_1$ and $K_2$. Therefore, by estimate \rf{estimate}, we have
$$\[\dbE \int_0^T  |\bar{X}(t)|_H^2 dt \]^\frac{1}{2} \les C \[\dbE \int_0^T  |\bar{\f}(t)|_H^2 dt \]^\frac{1}{2},$$
which implies the stability estimate.
\ef

\rm
The authors in \cite{Hamaguchi} gave the well-posedness of linear SVIEs by the method of continuation and discrete Gronwall inequality. In contrast, here we prefer to show the proof by induction directly. We also make comparisons with two important papers \cite{LY, Zhang} in the Appendix.

\br{comparison-SVIE}
To conclude this subsection, we point out that our singularity on the forward SVIEs can cover several important cases.
	Here are some of them.
	\begin{itemize}
		\item Let us consider the fractional singular kernel of the form $(t-s)^{-\a}, \a\in(0,\frac{1}{2})$, which is studied in Pr\"omel--Scheffels \cite{PS}, Hu et al. \cite{LHH}, to mention a few. As shown in Example \ref{frac}, this class of kernels can be covered by ours.
		\item For the convolution kernel  $K(t,s)\= h(t-s), \ (t,s)\in \D$ with $h(\cd) \in L^2(0,T;\dbR_+)$, it can be shown that $K(\cd)\in \sL^2(\D)$ by basic calculations. Observe that SVIEs with such kernels were applied in the Volterra Heston model (Abi Jaber \cite{Abi1}) which is a generalization of Heston stochastic volatility model in finance.
		\item SVIEs with completely monotone and convolution-type kernels can also be viewed as special case of ours, we refer to \cite[Appendix Lemma A.1]{Abi2} for more details.
	\end{itemize}
\er

\subsection{Several examples on singular stochastic Volterra systems}

In this subsection, we present several examples on singular backward stochastic Volterra integral equations. The first one shows a particular case of BSVIE whose singularity
improves the counterparts in (H3.1).

\begin{example}\label{frac}
Let us consider BSVIE (\ref{BSVIE}) where the coefficient $L_y(\cd,\cd)$ takes
the following case
$$\ba{ll}
\ns\ds L_y(t,s) = \frac{1}{(s-t)^{1-\b}}, \ \  (t,s) \in \D^*,\ \ \b\in(0,1).
\ea
$$
On the one hand, to keep $L_y \in L^2(\D^*)$ as in (H3.1),
we need $\b\in (\frac{1}{2},1)$. On the other hand, we can take advantage of this special framework and allow a big range for $\b$, i.e., $\b \in (0,1)$.
In fact, in this case, we can borrow the following Young's convolution inequality (by slightly modifying the proof of Corollary 2.2 in \cite{LY})
$$\(\int_0^T \lt| \int_t^T \frac{\f(s)}{(s-t)^{1-\beta}}ds \rt|_{\dbR}^pdt\)^\frac{1}{p}
	\les \( \frac{T^{1-r(1-\b)}}{1-r(1-\b)}\)^\frac{1}{r} \lt\|\f(\cd)\rt\|_{L^q(0,T;\dbR)},$$
where $\b \in (0,1)$, $p,q \ges 1$, $1\les r < \frac{1}{1-\b}$ satisfies $\frac{1}{p} +1 = \frac{1}{q} + \frac{1}{r}$, and $\f(\cd)\in L^q(0,T;\dbR)$.
Then choosing $p=q=2$, $r=1$, by the same procedure in \rf{1es} and  we estimate the $L_y$ term in \rf{0es} by the above inequality instead of H\"older ineqaulity, i.e.,
	$$ \dbE \int_S^T \( \int_t^T L_y(t,s)|y(s)|_Hds\)^2 dt = \dbE \int_S^T \(\int_t^T \frac{|y(s)|_H}{(s-t)^{1-\b}} ds\)^2 dt \les  \frac{(T-S)^{2\b}}{\b^2} \|y(\cd)\|^2_{L_\dbF^2(S,T;H)}.
	$$
	Similar arguments can be done to prove the contractivity of map $\Th$. Thus,
we can obtain the well-posedness of \rf{BSVIE} accordingly.

We emphasize that the above particular BSVIEs can be used in treating optimal control problems for forward SVIEs (\ref{fSVIE}) with $K_1$, $K_2$ in (H3.3) taking the following form (see e.g. \cite{LHH}),
	$$\ba{ll}
	\ns\ds K_1(t,s) = \frac{1}{(t-s)^{1-\a}}, \ \ K_2(t,s) = \frac{1}{(t-s)^{1-\g}}, \q (t,s) \in \D,\ \ \a,\ \g \in (0,\frac{1}{2}).
	\ea
	$$
For (\ref{fSVIE}), due to the above particular choice, we can allow $\a \in (0,1)$ by using Young's convolution inequality \cite[Corollary 2.2]{LY}.

\begin{comment}
{\color{red}

Put the following into the appendix or delete it.

More precisely, we first recall the Young¡¯s convolution inequality. Let $\b \in (0,1)$, $p,q \ges 1$ and $1\les r < \frac{1}{1-\b}$ be constants satisfying $\frac{1}{p} +1 = \frac{1}{q} + \frac{1}{r}$, then for any $\f(\cd)\in L^q(0,T;\dbR)$, it holds that
	$$\(\int_0^T \lt| \int_0^t \frac{\f(s)}{(t-s)^{1-\beta}}ds \rt|_{\dbR}^pdt\)^\frac{1}{p}
	\les \( \frac{T^{1-r(1-\b)}}{1-r(1-\b)}\)^\frac{1}{r} \lt\|\f(\cd)\rt\|_{L^q(0,T;\dbR)}.$$
	%
	Then, the well-posedness of \rf{fSVIE} can be obtained accordingy.
 }
\end{comment}

\end{example}

\medspace%

The following gives an example of linear BSVIE whose singularity can cover that in rough Heston model and fractional Brownian motion kernel properly.
\bex{Fractional Brownian motion Kernel}
For any constant $\cH\in (0,1)$, let
$$B^\cH(t)\= \int_0^t K_\cH(t,s) dB(s), \q t\in [0,T]. $$
Here $B(\cd)$ is a 1-dimensional standard Brownian motion and
$K_\cH(\cd,\cd)$ is defined by
$$K_\cH(t,s) \= \ell_1(t-s)^{\cH-\frac{1}{2}} + \ell_2 s^{\cH-\frac{1}{2}} F(\frac{t}{s}), \q (t,s)\in \D[0,T], $$
where
$\ds F(u)\= c_\cH(\frac{1}{2}-\cH)\int_1^u (r-1)^{\cH-\frac{3}{2}}(1-r^{\cH-\frac{1}{2}})dr,$ and %
$\ds c_\cH \= \(\frac{2\cH\G(\frac{3}{2}-\cH)}{\G(\cH+\frac{1}{2})\G(2-2\cH)}\)^\frac{1}{2}$, $\G(\cd)$ denotes the Gamma function, $\ell_1, \ell_2$ are two constants. In the following, we first look at two important special cases with different $\ell_1, \ell_2$.
\begin{itemize}
	\item[(i)] When $\ds \ell_1 = \frac{1}{\G(H+\frac{1}{2})}$ and $\ell_2=0$, $K_\cH^{(1)}(t,s) \= \ell_1(t-s)^{\cH-\frac{1}{2}}$ is called the Riemann-Liouville fractional kernel first introduced by L\'evy. At this moment, $K_\cH^{(1)}(\cd,\cd) \in \sL^2(\D) \subset L^2(\D)$. This kind of kernels is closely related to the rough Heston model, in for example \cite{BD,EER}.
	
	\item[(ii)] When $\ell_1=c_\cH$ and $\ell_2=1$, $\ds K_\cH^{(2)}(t,s) \= c_\cH (t-s)^{\cH-\frac{1}{2}} + s^{\cH-\frac{1}{2}} F(\frac{t}{s})$, $K_\cH^{(2)}(\cd,\cd)$ is the kernel of fractional Brownian motion $B^\cH(\cd)$ with Hurst parameter $\cH\in (0,1)$. By \cite[Theorem 3.2]{DU}, for any $\cH \in (0,1)$ we have
	\bel{es-fbk} K_\cH^{(2)}(t,s)\les C s^{-|H-\frac{1}{2}|} (t-s)^{-(\frac{1}{2}-H)_+},\q (t,s)\in \D, \ee
	where $C$ is a constant and $x_+ = max\{x,0\}$. It is easy to check
	$K_\cH^{(2)}(\cd,\cd) \in L^2(\D)$.
\end{itemize}

With the above kernels, let us look at the following linear BSVIE
\bel{bfbk-0}
Y(t) = \psi(t) + \int_t^T \[ K_\cH^{(2)} (s,t) M_1(t)^* Y(s) + K_\cH^{(1)} (s,t) M_2(t)^* Z(s,t) \]ds -\int_t^T Z(t,s) dW(s), \q t\in [0,T],
\ee
where $M_1(\cd), M_2(\cd)\in L^\infty(0,T;\cL(H;H))$.  By the estimate \rf{es-fbk} and Example \ref{d-s-k}, \rf{bfbk-0} admits a unique adapted M-solution according to Theorem \ref{43}. If we take
further particular form of $M_1(\cd)$, $M_2(\cd)$, then (\ref{bfbk-0}) can be used to tackle optimal control problem for the following forward stochastic Volterra equation:
\bel{fbk}X(t) =x + \int_0^t K_\cH^{(2)}(t,s)b(s,X(s))ds + \int_0^t K_\cH^{(1)}(t,s)\si(s,X(s))dW(s), \q t\in[0,T].\ee
Here $x\in H$ and $b,\;\si$ are two proper Lipschitz functions.
%
%$$ b(t,0) =0, \qq \si(t,0) = 0,$$
%$$|b(t,x)-b(t,y)|_H + |\si(t,x) - \si(t,y)|_{\cL_2^0} \les C|x-y|_H,$$ for all $t\in [0,T]$ and $x,y\in H$.
%
It is easy to see that \rf{fbk} admits a unique adapted solution $X(\cd)$ according to Theorem \ref{thm3}.
\ex

\ms
%{\color{red}The following example shows our BSVIE relates to the mild solution of the backward stochastic differential equation with Caputo derivative.}

\begin{comment}
\color{red} [Another definition of the mild solution in their paper]
$$\begin{aligned}
	y(t) =& E_\a(A(T-t)^\a) \xi + \int_t^T (s-t)^{\a-1} E_{\a,\a}(A(s-t)^\a) f(s,y(s),z(t,s)) ds \\
	%
	&\q + \int_t^T (s-t)^{\a-1} E_{\a,\a}(A(s-t)^\a) z(t,s) dw(s), \qq a.s., \ t\in [0,T],
\end{aligned}$$
where $E_\a(\cd)$ and $E_{\a,\a}(\cd)$ are one and two-parameter Mittag-Leffler functions,  \tb{is it the particular case in the infinite dimensional framework?} respectively, such that
%
$$\sup\{ \| E_\a(At^\a) \|, 0\les t\les T\} <\i, \q \sup\{ \| E_{\a,\a}(At^\a) \|, 0\les t\les T\} <\i.$$

%\no More details of Mittag-Leffler functions can be found in \cite{NA}.
\tb{ From \cite[Proposition 1.12, pp.12]{ZY}, we know that for $\a \in (0,1)$ and $t\in \dbR$, $E_{\a,\a}(t) > 0$.}
 Thus, define
$$Y(t)\=y(t), \q Z(t,s)\=(s-t)^{\a-1}E_{\a,\a}(A(s-t)^\a) z(t,s),$$
%
$$\left\{ \begin{aligned}
	& \psi(t) \= E_\a(A(T-t)^\a) \xi, \\
	%
	& g(t,s,Y(s),Z(t,s)) \= (s-t)^{\a-1} E_{\a,\a}(A(s-t)^\a) f\left(s,Y(s),(s-t)^{\a-1} E_{\a,\a}(A(s-t)^\a)^{-1} Z(t,s)\right).
\end{aligned}\right.$$
\end{comment}

\ms

Next, we give one example of infinite-dimensional linear BSVIEs motivated by optimal control problem for fractional stochastic evolution equations.

\bex{CFSEE}

To begin with, let us make some preparations. Suppose $q\in(\frac 1 2 ,1)$, and
$\cA$ is the infinitesimal generator of a uniformly bounded $C_0$-semigroup $\{Q(t)\}_{t\ges0}$ on $H$.
We define
$$S_q(t) \= \int_0^\i W_q(\th)Q(t^q\th)d\th, \q P_q(t) \= \int_0^\i q\th W_q(\th)Q(t^q\th)d\th,$$
where the Wright function
$$W_q(t) \=\frac{1}{\pi} \sum_{n=1}^{\infty} \frac{(-t)^{n}}{(n-1) !} \Gamma(n q) \sin (n \pi q), \quad t\ges 0,$$
satisfies
$$W_q(s) \ges 0, \q s\ges0,\qq \text{ and } \qq \int_0^\i W_q(s) ds = 1.$$
From \cite[Chapter 4, pp.136]{ZY}, we know that  $\{S_q(t)\}_{t>0}$ and $\{P_q(t)\}_{t>0}$ are uniformly bounded.

Consider the following linear BSVIE,
\bel{ad1}\begin{aligned}
Y(t) = \psi(t) + \int_t^T \[ (s-t)^{q-1} M_1(t)^* P_q(s-t)^* Y(s)  &+ (s-t)^{q-1} M_2(t)^* P_q(s-t)^* Z(s,t) \] ds \\
&\qq -\int_t^T Z(t,s)dW(s), \q t\in [0,T],
\end{aligned}\ee
where $M_i(\cd)\in L^\infty(0,T;
\cL(H;H))$.
Obviously, \rf{ad1} is an example of \rf{BSVIE}. Such kind of linear backward  Volterra integral equations plays an important role in the optimal control problems for the following fractional stochastic evolution equation with Caputo fractional derivative
\bel{CF} \left\{ \begin{aligned}
	&^C_0 D^q_t x(t) = \cA x(t) + f(t,x(t)) + g(t,x(t))\frac{dW(t)}{dt}, \q a.e.\; t\in [0,T], \\
	&x(0) =x_0,
\end{aligned}
\right.\ee
where $x_0 \in H$, $^C_0 D^q_t$ is the (left) Caputo fractional derivative of order $q\in ({\frac{1}{2}},1)$ defined by
$$^C_0D_t^q \f(t) \= \frac{1}{\G(1-q)}\frac{d}{dt} \int_0^t (t-s)^q \{\f(s) - \f(0)\}ds, \q t\in [0,T],$$
for appropriate function $\f$ (The integral above is understood in Bochner's sense).
In addition, the maps $f: [0,T] \times H \rightarrow H,\; g: [0,T] \times H \rightarrow \cL_2^0$ satisfy linear growth and Lipschitz conditions.

By the mild solution of \rf{CF}, we mean that the process $x(\cd) \in {L_\dbF^2(0,T;H)}$ which satisfies
\bel{CF-integral}
x(t) = S_q(t)x_0 + \int_0^t (t-s)^{q-1} P_q(t-s)f(s,x(s))ds + \int_0^t (t-s)^{q-1} P_q(t-s)g(s,x(s))dW(s), \q t\in [0,T],\ee
Obviously, the \rf{CF-integral} is  a
special case of \rf{fSVIE},
where $$\left\{ \begin{aligned}
	&\f(t)\= S_q(t)x_0,     \\
	&A(t,s,x(s)) \= (t-s)^{q-1} P_q(t-s)f(s,x(s)),        \\
	&B(t,s,x(s)) \=  (t-s)^{q-1} P_q(t-s)g(s,x(s)).        \\
\end{aligned}\right.$$

\no It's easy to see that the assumptions (H3.3) are satisfied and  \rf{CF} admits a unique mild solution according to Theorem \ref{thm3}.

We will give a more deep investigation in Section 4 for the optimal control problems for general  \rf{fSVIE}. As to the deterministic system in finite dimensional case, we refer to e.g. \cite{LY}.
\ex

To get more feelings about fractional partial differential equations, we give the following simple example.
\bex{}
Let $H=V\=L^2([0,\pi])$,
\bel{pde} \left\{ \begin{aligned}
	&^C_0D^q_t x(t,z) =  x_{zz}(t) + \sin x(t) + \cos x(t)\frac{dW(t)}{dt}, \q t\in [0,T], z\in [0,\pi], \\
	&x(t,0) = x(t,\pi) = 0, \q t\in [0,T],\\
	&x(0,z) =x_0(z), \q z\in [0,\pi],
\end{aligned}
\right.\ee
where $^C_0D^q_t$ is Caputo fractional derivative of order $q= \frac{3}{4}$, $x_0(\cd) \in H$.

Define an operator $\cA$ by $\cA x = x''$ with the domain
$$D(\cA) \= \{ x \in H \;\big|\; x,x' \text{ absolutely continuous, } x'' \in H, x(0) = x(\pi) = 0.\}$$
Then $\cA$ generates a $C_0$-semigroup $\{Q(t)\}_{t\ges 0}$ which is compact, analytic and self-adjoint. Moreover, $Q(\cd)$ is a uniformly stable semigroup and $\|Q(t)\| \les 1.$
In fact, $\cA$ has a discrete spectrum, the eigenvalues are $\{-n^2\}_{n\in \dbN}$, and the corresponding orthogonal eigenvectors are given by $e_n(z) = \sqrt{\frac{2}{\pi}} \sin(nz)$. And for each $v\in H$,
$$ \cA v= \sum_{n=1}^\i n^2 \lan v, e_n \ran e_n, \q
Q(t)v =  \sum_{n=1}^\i e^{-n^2t}\lan v,e_n \ran e_n.$$
Let $f(t,x)\= \sin x, \; g(t,x) \= \cos x$, then the assumptions on $f,g$ in Example \ref{CFSEE} are fulfilled. %Therefore, \rf{pde} admits a unique mild solution.
\ex

\ms

\section{Two classes of optimal control problems}

As applications, in this section we study the optimal control problems for controlled stochastic Volterra integral equation and controlled stochastic delay evolution equation, respectively. We present the maximum principles of optimal control with convex control regions.

\subsection{Maximum principle for controlled stochastic Volterra integral equations}

In this subsection, motivated by the examples of fractional stochastic evolution equations and stochastic evolutionary integral equations, we discuss an optimal control problem for a stochastic Volterra integral equation with a Lagrange type cost functional. More precisely, consider the following state equation:
\bel{state-0}
X(t) = \f(t) + \int_0^t b(t,s,X(s),u(s)) ds + \int_0^t \si(t,s,X(s),u(s)) dW(s), \q t\in [0,T],
\ee
with the cost functional
\bel{cost}
J(u(\cd)) = \dbE \int_0^T g(t,X(t),u(t)) dt.
\ee
Here $u(\cd)$ is the {\it control process} valued in $U$, with $U$ a non-empty convex set of a separable metric space $\cO$, and $X(\cd)$ is the corresponding {\it state process} valued in $H$. In the above, $\f: [0,T] \times \O \mapsto H$ is called the {\it free term}. $b:\D \times H \times U \times \O \mapsto H$ and $\si: \D \times H \times U \times \O \mapsto \cL_2^0$ are called the {\it drift and diffusion } of the controlled SVIEs, respectively, $g: [0,T] \times H \times U \times \O \mapsto \dbR$ is called the {\it running cost}.
We now introduce the following assumptions.\\

{\bf (H4.1)} Suppose $b, \si $, $g$ are continuously differentiable	in the variables $(x,u)$, and the Fr\'echet derivatives $g_x$ (over $H$) and $g_u$ (over $U$) are bounded. Moreover,
$$ \ba{ll}
\ns\ds |b(t,s,0,u)|_H + |\si(t,s,0,u)|_{\cL_2^0} \les C,   \\
\ns\ds |b_x(t,s,x,u)|_{L(H)} + |b_u(t,s,x,u)|_{L(U;H)} \les K_1(t,s),\\
\ns\ds |\si_x(t,s,x,u)|_{L(H;\cL_2^0)} + |\si_u(t,s,x,u)|_{L(U;\cL_2^0)} \les K_2(t,s),
\ea$$
for any $(t,s)\in \D,\; x\in H, \; u \in U,$ and $K_1(\cd,\cd) \in L^2(\D)$, $K_2(\cd,\cd) \in \sL^2(\D)$.\\

Let $\sU[0,T] \= L_\dbF^2(0,T;U)$
be the set of admissible controls. Under {\bf (H4.1)}, by Theorem \ref{thm3}, for any $\f(\cd) \in L_\dbF^2(0,T;H)$ and $u(\cd) \in \sU[0,T]$, \rf{state-0} admits a unique solution $X(\cd) \in L_\dbF^2(0,T;H)$. Then the cost functional $J(u(\cd))$ is well-defined and the optimal control problem can be stated as follows:\\

{\bf Problem (C).} Find a $\bar{u}(\cd)\in \sU[0,T]$ such that
\bel{problem}
J(\bar{u}(\cd)) = \inf_{u(\cd) \in \sU[0,T]} J(u(\cd)).
\ee
In this case, $\bar{u}(\cd)$ is called an {\it optimal control}, the corresponding state process $\bar{X}(\cd)$ and $(\bar{X}(\cd), \bar{u}(\cd))$ are called the {\it optimal state } and {\it optimal pair}, respectively.

\bt{MP} Let $(\bar{X}(\cd), \bar{u}(\cd))$ be an optimal pair. Then the following BSVIE admits a unique adapted $M$-solution  $(Y(\cd),Z(\cd,\cd))$,
\bel{adjoint} \begin{aligned}
Y(t) &= g_x(t,\bar{X}(t), \bar{u}(t)) + \int_t^T \Big( b_x(s,t,\bar{X}(t), \bar{u}(t))^*Y(s) \\
&\qq\qq \qq\qq + \si_x(s,t,\bar{X}(t), \bar{u}(t))^*Z(s,t) \Big)ds - \int_t^T Z(t,s)dW(s), \q t\in [0,T],
\end{aligned} \ee
such that
\bel{conclusion}\begin{aligned}
\Big\lan g_u(t,\bar{X}(t), \bar{u}(t)) +\dbE\Big[ &\int_t^T b_u(s,t,\bar{X}(t), \bar{u}(t))^*Y(s)ds \\
&+ \int_t^T \si_u(s,t,\bar{X}(t), \bar{u}(t))^*Z(s,t)ds\Big|\cF_t \Big],
 u - \bar{u}(t)\Big\ran_U \ges 0, \q \forall u\in U, \; t\in[0,T], \; \as
\end{aligned} \ee
\et

\pf To begin with, we introduce the following abbreviations:
$$\lt\{ \begin{aligned}
	&b_x(t,s)\= b_x(t,s,\bar{X}(s), \bar{u}(s)), \qq b_u(t,s)\= b_u(t,s,\bar{X}(s),\bar{u}(s)),\\
	&\si_x(t,s)\= \si_x(t,s,\bar{X}(s), \bar{u}(s)), \qq \si_u(t,s)\= \si_u(t,s,\bar{X}(s),\bar{u}(s)),\\
    &g_x(t) \= g_x(t,\bar{X}(t), \bar{u}(t)), \qq g_u(t) \= g_u(t,\bar{X}(t), \bar{u}(t)).
\end{aligned}\rt.$$
Let $( \bar{X}(\cd), \bar{u}(\cd) )$ be an optimal pair of Problem (C). For any $v(\cd) \in \sU[0,T]$ and any $\e \in \dbR$, denote
$$u^\e(\cd) \= \bar{u}(\cd) + \e [v(\cd) - \bar{u}(\cd)] \in \sU[0,T].$$
Let $X^\e(\cd)$ be the solution of \rf{state-0} with $u(\cd)$ replaced by $u^\e(\cd)$. Then

$$  \begin{aligned}
	X^\e(t) - \bar{X}(t) =&\int_0^t \(  \tilde{b}_x(t,s)[X^\e(s) - \bar{X}(s)]  +  \tilde{b}_u(t,s) \e [v(\cd) - \bar{u}(\cd)] \) ds \\
	& \qq + \int_0^t \(  \tilde{\si}_x(t,s)[ X^\e(s) - \bar{X}(s)]  +   \tilde{\si}_u(t,s) \e [v(\cd) - \bar{u}(\cd)] \) dW(s),
\end{aligned} $$
where
$$ \begin{aligned}
	&\ti{b}_x(t,s) \= \int_0^1 b_x(t,s,\bar{X}(s) + \a [X^\e(s) - \bar{X}(s)], u^\e(s)) d\a,\\
	&\ti{b}_u(t,s) \= \int_0^1 b_u(t,s,\bar{X}(s),\bar{u}(s)+\th \e v(s)) d\th, \\
	&\ti{\si}_x(t,s) \= \int_0^1 \si_x(t,s,\bar{X}(s) + \a [X^\e(s) - \bar{X}(s)], u^\e(s)) d\a,\\
	&\ti{\si}_u(t,s) \= \int_0^1 \si_u(t,s,\bar{X}(s),\bar{u}(s)+\th \e v(s)) d\th.
\end{aligned} $$
By the estimate \rf{estimate} in Theorem \ref{thm3}, we have
$$\begin{aligned}
	\lt\| X^\e(\cd) - \bar{X}(\cd) \rt\|^2_{L^2_\dbF(0,T;H)} &\les C \dbE \int_0^T \lt| \int_0^t  \tilde{b}_u(t,s) \e [v(\cd) - \bar{u}(\cd)]  ds + \int_0^t \tilde{\si}_u(t,s) \e [v(\cd) - \bar{u}(\cd)] dW(s) \rt|_H^2 dt \\
	&\les C\e^2.
\end{aligned}$$
Hence,
$$\lim_{\e \rightarrow 0} \lt\| X^\e(\cd) - \bar{X}(\cd) \rt\|_{L^2_\dbF(0,T;H)} = O(\e).$$

\no Define
$$ X^\e_1(t) \= \frac{X^\e(t) - \bar{X}(t)}{\e}, \qq t\in [0,T].$$
Then $ X^\e_1(\cd) \rightarrow X_1(\cd)$ in $L_\dbF^2(0,T;H)$ with $X_1(\cd)$ satisfying the following equation:
\bel{variational} \begin{aligned}
  X_1(t) &= \int_0^t \{  b_x(t,s) X_1(s)  +  b_u(t,s)[v(s) - \bar{u}(s)] \}ds \\
  &
  \qq\qq\qq + \int_0^t \{ \si_x(t,s) X_1(s)  +  \si_u(t,s) [v(s) - \bar{u}(s)] \}dW(s)  \\
  &\equiv \bar{\f}(t) + \int_0^t b_x(t,s) X_1(s)ds + \int_0^t \si_x(t,s) X_1(s)dW(s), \qq t\in [0,T],
\end{aligned}\ee
where
$$\bar{\f}(t) \=\int_0^t b_u(t,s)[v(s) - \bar{u}(s)] ds + \int_0^t \si_u(t,s)[v(s) - \bar{u}(s)] dW(s), \q t\in[0,T].$$\\

\no By the optimality of $( \bar{X}(\cd), \bar{u}(\cd) )$, we have
$$\begin{aligned}
	0 &\les \frac{J(u^\e(\cd))- J(\bar{u}(\cd))} {\e} \\
	&= \dbE \int_0^T \[\llan g_x(t,\bar{X}(t)+\a[X^\e(t) - \bar{X}(t)], u^\e(t)), \frac{X^\eps(t) - \bar{X}(t)}{\e} \rran_H \\
	&\qq \qq \qq \qq + \lan g_u(t,\bar{X}(t),\bar{u}(t) +\th \e [v(t) - \bar{u}(t)] ), v(t) - \bar{u}(t) \ran_U \] dt \\
	&\rightarrow \dbE \int_0^T \[\lan g_x(t), X_1(t) \ran_H + \lan g_u(t), v(t) - \bar{u}(t) \ran_U \] dt \\
	&\equiv I_1 + I_2,
\end{aligned}$$
where $\a,\th \in (0,1)$ and $\ds I_1 \=\dbE \int_0^T \lan g_x(t), X_1(t) \ran_H dt $ and $\ds I_2 \= \dbE \int_0^T \lan g_u(t), v(t) - \bar{u}(t) \ran_U dt $.

\no To get rid of the term involving $X_1(\cd)$, we have
$$\begin{aligned}
	&\dbE \int_0^T \lan \bar{\f}(t), Y(t) \ran_H dt\\
	=& \dbE \int_0^T\[ \llan X_1(t) - \int_0^t b_x(t,s)X_1(s)ds - \int_0^t \si_x(t,s)X_1(s)dW(s), Y(t) \rran_H \] dt\\
	=&\dbE \int_0^T\[ \llan X_1(t),Y(t)\rran_H - \int_0^t \llan b_x(t,s)X_1(s), Y(t)\rran_H ds - \int_0^t \llan \si_x(t,s)X_1(s), Z(t,s)\rran_{\cL_2^0} ds\]dt \\
	=& \dbE \int_0^T \llan X_1(t),Y(t)\rran_H dt - \dbE \int_0^T \llan X_1(t), \int_t^T b_x(s,t)^*Y(s)ds + \int_t^T \si_x(s,t)^*Z(s,t)ds \rran_H dt \\
	=& \dbE \int_0^T \llan X_1(t), Y(t) - \int_t^T b_x(s,t)^*Y(s)ds - \int_t^T \si_x(s,t)^*Z(s,t)ds \rran_H dt \\
	=& \dbE \int_0^T \llan X_1(t), g_x(t) - \int_t^T Z(t,s)dW(s) \rran_H dt\\
	=& \dbE \int_0^T \llan X_1(t), g_x(t)\rran_H dt.
\end{aligned} $$
{The above equality is the so-called duality principle between \rf{adjoint} and \rf{variational}.}
The first equality holds because of \rf{variational}. By the definitions of $M$-solution and adjoint operator, we get the second and third equalities via Fubini theorem. The fifth equality holds due to \rf{adjoint}. The properties of stochastic integrals imply the last equality. By the definition of $\bar{\f}(\cd)$, we derive that
$$ \begin{aligned}
	I_1 &\equiv \dbE \int_0^T \lan  g_x(t),X_1(t)\ran_H dt
	= \dbE \int_0^T \lan \bar{\f}(t), Y(t) \ran_H dt \\
	&= \dbE \int_0^T \lt\lan \int_t^T b_u(s,t)^*Y(s)ds + \int_t^T \si_u(s,t)^*Z(s,t)ds, v(t) - \bar{u}(t) \rt\ran_U dt.
\end{aligned}$$
As a consequence,
$$\begin{aligned}
	0 &\les I_1 +I_2 \\
	=& \dbE \int_0^T \lt\lan g_u(t) + \int_t^T b_u(s,t)^*Y(s)ds + \int_t^T \si_u(s,t)^*Z(s,t)ds, v(t) - \bar{u}(t) \rt\ran_U dt \\
	%
	%=& \dbE \int_0^T \lt\lan g_u(t) + \int_t^T b_u(s,t)^*Y(s)ds + \int_t^T \si_u(s,t)^*Z(s,t)ds - \int_t^T Z(t,s)dW(s), v(t) - \bar{u}(t) \rt\ran dt\\
	%
    %=& \dbE \int_0^T \lt\lan Y(t), v(t) - \bar{u}(t) \rt\ran dt.
\end{aligned}$$
which implies \rf{conclusion} since $v(\cd)$ is arbitrary.
\ef
 When the cost functional \rf{cost} depends on both the running cost and  terminal cost, we need the solution $X(\cd)$ of \rf{state-0} to be (left) continous at $t=T$ and make similar arguments. We leave it to the interested readers.

\subsection{Maximum principle for controlled stochastic delay evolution equations}

 In this part, we study the optimal control problems for controlled stochastic delay evolution equations (SDEEs, for short) by the theory of BSVIEs developed in Section 3.

Given a constant time delay parameter $\d\in(0,T)$ and a constant $\l\in\dbR$, $U$ is a nonempty convex subset of a separable Hilbert space $\fU$. In this subsection, we consider the controlled stochastic delay evolution equation of the following form:
\vskip-7mm
\bel{state-}\left\{\ba{ll}
\ds d x(t)= \cA x(t)dt + b\(t,x(t),x(t-\d),\int_{-\d}^0e^{\l\th}x(t+\th)d\th,u(t),u(t-\d)\)dt\\
\ns\ds\q \qq+\si\(t,x(t),x(t\1n-\1n\d),\2n\int_{-\d}^0\3ne^{\l\th}x(t+\th)d\th,\1nu(t), \1nu(t\1n-\1n\d)\)dW(t),\ t\in[0,T],\\
\ns\ds x(t)=\xi(t),\ u(t)=\eta(t),\ t\in[-\d,0],
\ea\right.\ee
\vskip-1mm
\no where $x(\cd)\in H$ is state and $u(\cd)\in U$ is control, $\cA$ is the infinitesimal generator of a $C_0$-semigroup $\{S(t)\}_{t\ges 0}$ on $H$,  %Suppose that $(\Omega,\mathcal{F},\dbF,\mathbb{P})$ is a complete filtered probability space and the filtration $\dbF=\{\mathcal{F}_t\}_{t\geq 0}$ is generated by a $d$-dimensional standard Brownian motion $\{W(t)\}_{t\geq0}$.
$b:[0,T] \times H \times H \times H \times U \times U \rightarrow H $, $\si: [0,T] \times H \times H \times H \times U \times U \rightarrow \cL_2^0 $ are given maps.
Deterministic continuous function $\xi(\cd)$ and square integrable function $\eta(\cd)$ are the initial trajectories of the state and the control, respectively. We associate \rf{state-} with the following cost functional
%
%\vskip-3mm
%
\bel{cost-}\ba{ll}
\ns\ds J(u(\cd))=\dbE \int_0^T l\(t,x(t),x(t-\d),\int_{-\d}^0e^{\l\th}x(t+\th)d\th,u(t),u(t-\d)\)dt \\
\ns\ds \qq\qq\qq+h\(x(T),x(T-\d),\int_{-\d}^0e^{\l\th}x(T+\th)d\th\)\bigg], \ea\ee
\vskip-1mm
\no where $l: [0,T] \times H \times H \times H \times U \times U \rightarrow \dbR $ and $h: H \times H \times H \rightarrow \dbR $ are given.
%$b:[0,T]\times\Om\times\dbR^n\times\dbR^n\times\dbR^n\times U\times U \rightarrow\dbR^n$, $\sigma:[0,T]\times\Om\times\dbR^n\times\dbR^n\times\dbR^n\times U\times U \rightarrow\dbR^{n\times d}$, $l:[0,T]\times\Om\times\dbR^n\times\dbR^n\times\dbR^n\times U\times U \rightarrow\dbR$ and $h:\Om\times\dbR^n\times\dbR^n\times\dbR^n\rightarrow\dbR$ are given mappings.
%
Define the admissible control set as follows:
%
%\vskip-7mm
%
$$\ba{ll}
\ns\ds \mathcal{U}_{ad}\=\Big\{u(\cdot):[-\d,T]\rightarrow U \big|\; u(\cdot)\mbox{ is a }U\mbox{-valued, square-integrable},\ \dbF\mbox{-adapted}\\
\ns\ds  \qq\qq\qq\qq\qq\qq\q \mbox{process and }u(t)=\eta(t),\ t\in[-\d,0]\Big\}.
\ea$$
\vskip-3mm

The optimal control problem is stated as follows:\\

\textbf{Problem (P).} Find a control $u^*(\cdot)$ over $\mathcal{U}_{ad}$ such that (\ref{state-}) is satisfied and (\ref{cost-}) is minimized, i.e.,
$$\ba{ll}
J(u^*(\cd))=\inf\limits_{u(\cd)\in\,\mathcal{U}_{ad}}J(u(\cd)).
\ea$$
\vskip 2mm

In order to rewrite the state equation \rf{state-} in a more concise form, denote
\bel{yzm}\ba{ll}
\ns\ds y(t)\=x(t-\d),\q z(t)\=\int_{-\d}^0e^{\l\th}x(t+\th)d\th,\q \m(t)\=u(t-\d).
\ea\ee
\vskip-3mm
\no Here $y(\cd)$ and $z(\cd)$ are considered to be the  point-wise delayed state and moving-average delayed state of $x(\cd)$, respectively. And $\mu(\cd)$ is called the point-wise delayed control. Then, the state equation \rf{state-} can be written as follows:
\vskip-4mm
\bel{state}\left\{\ba{ll}
\ns\ds d x(t)= \cA x(t)dt + b(t,x(t),y(t),z(t),u(t),\m(t))dt\\
\ns\ds\qq \qq+\si(t,x(t),y(t),z(t),u(t),\m(t))dW(t),\ \ t\in[0,T],\\
\ns\ds x(t)=\xi(t),\ u(t)=\eta(t),\ t\in[-\d,0].
\ea\right.\ee
\vskip-1mm
\no And the cost functional \rf{cost-} becomes
\vskip-4mm
\bel{cost}\ba{ll}
\ns\ds J(u(\cd))=\dbE\bigg[ \int_0^T l(t,x(t),y(t),z(t),u(t),\mu(t))dt +h(x(T),y(T),z(T))\bigg].
\ea\ee
We give the following assumptions.
\vskip 4mm
{\bf (H4.2)} The maps $(x,y,z) \mapsto b=b(t,x,y,z,u,\mu),\; \si= \si(t,x,y,z,u,\mu),\;l=l(t,x,y,z,u,\mu),\; h=h(x,y,z)$ are G\^ateaux differentiable in $(x,y,z)$ with continuous and bounded derivatives. And there exists a constant $C$ such that
$$ |b(t,0,0,0,u,\mu)|_H+ |\si(t,0,0,0,u,\mu)|_{\cL_2^0} \les C, \q \forall u,\mu\in U,\; t\ges 0.$$

Under (H4.2), the state equation \rf{state-} (or \rf{state}) admits a unique solution (see e.g. \cite{ZL}), so the cost functional \rf{cost-} (or \rf{cost}) is well-defined.
For notational simplicity, for $t \in [0,T]$, we denote $\Th(t) \= (x^*(t), y^*(t),z^*(t),u^*(t),\mu^*(t))$, $\Pi(t) \= (x^*(t), y^*(t),z^*(t))$. For any $v(\cd) \in \cU_{ad}$, $0 \les \e \les 1$, denote  $u^\e(t) \= u^*(t) + \e (v(t) - u^*(t)) \in \cU_{ad} $ and $\mu^\e(t) \= u^\e(t-\d)$. In addition, for $f\=b, \si, l$,
$$\begin{aligned}
	& f_x(t)\= f_x(t,\Th(t)),\q f_y(t)\= f_y(t,\Th(t)), \q f_z(t)\= f_z(t,\Th(t)), \\
	&h_x(T)\= h_x(\Pi(T)), \q h_y(T)\= h_y(\Pi(T)), \q h_z(T)\= h_z(\Pi(T)), \\
	& \D f(t) \= f (t,x^*(t), y^*(t),z^*(t),u^\e(t),\mu^\e(t)) - f(t,\Th(t)).
\end{aligned}$$

We introduce the variational equation (\cite{LZ}):
\vskip-5mm
\bel{variational equation-1}\left\{\ba{ll}
\ns\ds d x_1(t)=\Big[\cA x_1(t) +  b_x(t)x_1(t)+b_{y}(t)y_1(t)+b_z(t)z_1(t)+\D b(t)\Big]dt\\
\ns\ds \qq\q +\Big[\si_x(t)x_1(t)\1n+\1n\si_{y}(t)y_1(t) \1n+\1n\si_z(t)z_1(t)\1n+\1n\D \si(t)\Big]\1ndW(t), \q t\1n\in\1n[0,T],\\
\ns\ds x_1(t)=0,\q t\in[-\d,0],
\ea\right.\ee
where $y_1(\cd)$, $z_1(\cd)$ are defined similar to \rf{yzm}, i.e.
$$y_1(t) \= x_1(t-\d), \q z_1(t)\= \int_{-\d}^0 e^{\l\th} x_1(t+\th) d\th.$$

Next, we introduce the following variational inequality.

\bl{lemma variational inequality}
\no Let Assumption {\rm {\bf (H4.2)} } hold. Suppose $(x^*(\cdot),u^*(\cdot))$ is an optimal pair, $x^\varepsilon(\cd)$ is the trajectory corresponding to $u^\varepsilon(\cdot)$. Then, the following variational inequality holds:
\bel{variational-inequality}\ba{ll}
\ns\ds  J(u^\e(\cd))\1n-\1nJ(u^*(\cd))\1n= \dbE \Big[\lan h_x(T),x_1(T)\ran_H +\lan h_y(T),y_1(T)\ran_H+ \lan h_z(T),z_1(T)\ran_H \big] \\
\ns \ds \qq\qq\qq\qq+
\dbE\int_0^T \Big[\D l(t)+ \lan l_x(t),x_1(t) \ran_H +\lan l_y(t),y_1(t)\ran_H  +\lan l_z(t),z_1(t)\ran_H\Big]dt
+o(\e).
\ea \ee \el
\vskip-1mm

The proof is similar to  \cite[Theorem 12.4]{LZ}, so we omit the details here.
Inspired by \cite{MSWZ}, we rewrite the variational equation \rf{variational equation-1} into a linear SVIE. To this end, for any $t\in [0,T]$ we define
\vskip-4mm
$$\ba{ll}
\ns\ds X(t)\=\Bigg[\begin{array}{ccccc}
	x_1(t) \\
	y_1(t){\bf 1}_{(\d,\i)}(t) \\
	z_1(t)
\end{array}\Bigg] {\in \widetilde H \= H\times H \times H}, \ea $$

\no and for any $(t,s)\in \D$,%for $j=1,\cds,d$,
\vskip-7mm
$$
\ba{ll}
\ns\ds A(t,s)\=\Bigg[\begin{array}{ccccc}
	S(t-s)b_x(s) & S(t-s)b_y(s) & S(t-s)b_z(s) \\
	{\bf 1}_{(\d,\infty)}(t-s) S(t-s)b_x(s) & {\bf 1}_{(\d,\infty)}(t-s) S(t-s)b_y(s) & {\bf 1}_{(\d,\infty)}(t-s) S(t-s)b_z(s)\\
	{\cI} & -e^{-\l\d}\cI & -\l \cI
\end{array}\Bigg],\\
\ns\ds C(t,s)\=\Bigg[\begin{array}{ccccc}
	S(t-s)\si_x(s) & S(t-s)\si_y(s) & S(t-s)\si_z(s) \\
	{\bf 1}_{(\d,\infty)}(t-s) S(t-s)\si_{x}(s) & {\bf 1}_{(\d,\infty)}(t-s) S(t-s)\si_y(s) & {\bf 1}_{(\d,\infty)}(t-s) S(t-s)\si_z(s) \\
	\cO & \cO & \cO
\end{array}\Bigg],\\
B(t,s)\=\Bigg[\begin{array}{ccccc}
	S(t-s)\D b(s) \\
	{\bf 1}_{(\d,\infty)}(t-s) S(t-s)\D b (s) \\
	\cO
\end{array}\Bigg],\q
D(t,s)\=\Bigg[\begin{array}{ccccc}
	S(t-s)\D \si(s) \\
	{\bf 1}_{(\d,\infty)}(t-s) S(t-s)\D \si (s) \\
	\cO
\end{array}\Bigg],\ea$$
where $\cI,\cO$ denote the identity operator and zero operator on $H$, respectively.

Therefore, the variational equation \rf{variational equation-1} can be rewritten as
\bel{variational X} \begin{aligned}
	X(t) =& \int_0^t [A(t,s)X(s) + B(t,s)]ds + \int_0^t [C(t,s)X(s) + D(t,s)]dW(s),  \\
	\equiv&\ \f(t) + \int_0^t A(t,s)X(s)ds + \int_0^t C(t,s)X(s) dW(s),\q t\in[0,T],
\end{aligned}
\ee
where $$\f(t) \= \int_0^t B(t,s)ds + \int_0^t D(t,s) dW(s).$$
And the variational inequality \rf{variational-inequality} can be represented as
\bel{variational-inequality-X}
J(u^\e(\cd))-J(u^*(\cd))=  \dbE\int_0^T  \Big[\lan \bar L(t),X(t)\ran{_{\widetilde H}}  +\D l(t)\Big]dt +\dbE \lan \bar H, X(T) \ran{_{\widetilde H}}
+o(\e),
\ee

\no where $$\ba{ll}
%
%\ns\ds \bar H=\[\begin{array}{ccccc}
	%
	%h_x(T) & h_y(T) & h_z(T)
	%
%\end{array}\],\q
\bar L(t)\=\lt[\begin{array}{ccccc}
	l_x(t)\\ l_y(t) \\ l_z(t)
\end{array}\rt] \in \widetilde{H},\q \bar H\=\lt[\begin{array}{ccccc}
h_x(T)\\ h_y(T) \\ h_z(T)
\end{array}\rt] \in \widetilde{H}. \ea $$

Notice that \rf{variational X} is a linear SVIE, we introduce the following adjoint equation, %which is a BSVIE other than an anticipated BSDE,
\bel{adjoint-equ}\lt\{\begin{aligned}
& \bar \eta(t) = \bar H - \int_t^T \zeta(s) dW(s), \q t\in [0,T], \\
&Y(t) = \bar L(t) + A(T,t)^*\bar H + C(T,t)^*\zeta(t) + \int_t^T [A(s,t)^* Y(s) \\
& \qq\qq\qq\qq\qq + C(s,t)^* Z(s,t)] ds - \int_t^T Z(t,s) dW(s), \q t\in [0,T], \\
& Y(t) = \dbE Y(t) + \int_0^t Z(t,s)dW(s), \q
t\in [0,T],
\end{aligned}\rt.\ee
where for $\forall t\in [0,T],$
$$ \zeta(t)\= \lt[ \begin{aligned}
	\zeta^0(t) \\ \zeta^1(t) \\ \zeta^2(t)
\end{aligned}\rt] \in \widetilde{\cL_2^0} \= \cL_2^0 \times \cL_2^0 \times \cL_2^0, \q
Y(t)\= \lt[ \begin{aligned}
	Y^0(t) \\ Y^1(t) \\ Y^2(t)
\end{aligned}\rt] \in \widetilde{H}, \q Z(t,s) \= \lt[ \begin{aligned} Z^0(t,s) \\ Z^1(t,s)\\ Z^2(t,s)\end{aligned}\rt] \in \widetilde{\cL_2^0}. $$

\no By the dual principle similar to \cite[Theorem 5.1]{Y 2008}, we have
$$\dbE \int_0^T \lan \bar L(t), X(t) \ran{_{\widetilde H}} dt +\dbE \lan \bar H, X(T) \ran{_{\widetilde H}} = \dbE \int_0^T \lan \f(t), Y(t) \ran{_{\widetilde H}} dt +\dbE \lan \bar H, \f(T) \ran{_{\widetilde H}}.$$

\no Thus,
$$ \begin{aligned}
	& \dbE \int_0^T \lan \f(t), Y(t) \ran{_{\widetilde H}} dt \\
	=\ & \dbE \int_0^T \[ \llan \int_0^t B(t,s)ds, Y(t) \rran{_{\widetilde H}} + \llan \int_0^t D(t,s)dW(s), \dbE Y(t) + \int_0^t Z(t,s)dW(s) \rran{_{\widetilde H}}\] dt \\
	=\ & \dbE\int_0^T \[ \llan \int_0^t B(t,s)ds, Y(t) \rran{_{\widetilde H}} dt + \int_0^t \lan D(t,s),Z(t,s) \ran_{\widetilde{\cL_2^0}} ds \]dt \\
	=\ & \dbE\int_0^T \[ \int_{s}^{T} \lan B(t,s), Y(t) \ran{_{\widetilde H}} dt + \int_s^T \lan D(t,s), Z(t,s) \ran_{\widetilde{\cL_2^0}} dt \] ds \\
	=\ & \dbE\int_0^T \Big\{ \int_s^T \[ \lan S(t-s)\D b(s), Y^0(t)\ran_H + \lan S(t-s){\bf 1}_{(\d,\i)}(t-s) \D b(s),Y^1(t)\ran_H \] dt \\
	&\qq + \int_s^T \[ \lan S(t-s)\D \si(s), Z^0(t,s)\ran_{\cL_2^0} + \lan S(t-s){\bf 1}_{(\d,\i)}(t-s) \D \si(s),Z^1(t,s)\ran_{\cL_2^0} \] dt \Big\} \\
	=\ & \dbE\int_0^T \Big\{ \llan \D b(s), \int_s^T S(t-s)^*Y^0(t)dt + \int_{s+\d}^T S(t-s)^*Y^1(t)dt \cd {\bf 1}_{[0,T-\d)}(s) \rran_H \\
	&\qq + \llan \D \si(s), \int_s^T S(t-s)^*Z^0(t,s)dt + \int_{s+\d}^T S(t-s)^*Z^1(t,s)dt \cd {\bf 1}_{[0,T-\d)}(s) \rran_{\cL_2^0} \Big\} ds. \\
	%
	%\equiv\ & \dbE\int_0^T \[ \lan \D b(s), p(s) \ran_H  + \lan \D \si(s), q(s) \ran_{\cL_2^0} \] ds.
\end{aligned} $$

\no Similarly,
$$\begin{aligned}
\dbE \lan \bar H, \f(T) \ran{_{\widetilde H}} &= \dbE \int_0^T \[ \lan \bar H, B(T,s) \ran{_{\widetilde H}} + \lan \zeta(s), D(T,s) \ran_{\widetilde{\cL_2^0}} \] ds \\
&=\dbE \int_0^T \[ \llan \D b(s), S(T-s)^*h_x(T) + S(T-s)^*h_y(T)\cd {\bf 1}_{[0,T-\d)}(s) \rran_H \\
&\qq\qq\qq\qq+ \llan \D \si(s), S(T-s)^*\zeta^0(s) + S(T-s)^*\zeta^1(s)\cd {\bf 1}_{[0,T-\d)}(s)\rran_{\cL_2^0}\]ds.
\end{aligned}$$

Denote
$$\left\{\begin{aligned}
	&p(s) \=S(T-s)^*h_x(T) + S(T-s)^*h_y(T)\cd {\bf 1}_{[0,T-\d)}(s) \\
	 &\qq\qq\qq\qq\qq\qq + \int_s^T S(t-s)^*Y^0(t)dt  + \int_{s+\d}^T S(t-s)^*Y^1(t)dt \cd {\bf 1}_{[0,T-\d)}(s),\\
	&q(s) \=S(T-s)^*\zeta^0(s) + S(T-s)^*\zeta^1(s)\cd {\bf 1}_{[0,T-\d)}(s) \\
	&\qq\qq\qq\qq\qq+ \int_s^T S(t-s)^*Z^0(t,s)dt + \int_{s+\d}^T S(t-s)^*Z^1(t,s)dt \cd {\bf 1}_{[0,T-\d)}(s), \q s\in [0,T].
\end{aligned}\rt.$$

\no Then,
$$J(u^\e(\cd))-J(u^*(\cd))=  \dbE\int_0^T  \Big[\D l(s) + \lan \D b(s), p(s) \ran_H  + \lan \D \si(s), q(s) \ran_{\cL_2^0} \Big]ds+o(\e).$$

Next, define a Hamilton function $G: [0,T] \times H \times H \times H \times U \times U \times H \times \cL_2^0\rightarrow \dbR $ as
$$G(t,x,y,z,u,\mu,p,q) \= l(t,x,y,z,u,\mu) + \lan p, b(t,x,y,z,u,\mu) \ran_H + \lan q, \si(t,x,y,z,u,\mu) \ran_{\cL_2^0}.$$

\bt{} Let (H4.2) hold. If $u^*(\cd) \in \cU_{ad}$ is an optimal control, then the following maximum condition holds:
\bel{mp}\dbE \int_0^T \[ \lan G_u(t,\Th(t),p(t),q(t)),  v(t) - u^*(t) \ran_U + \lan G_\mu(t,\Th(t),p(t),q(t)),  v(t-\d) - u^*(t-\d) \ran_U \] dt \ges 0, \q \forall v(\cd) \in \cU_{ad}.\ee
\et

\pf By \rf{variational-inequality}, when $\e \to 0$,
$$\begin{aligned}
\lim_{\e \rightarrow 0} \frac{J(u^\e(\cd))-J(u^*(\cd))}{\e} =& \lim_{\e \to 0} \frac{1}{\e} \dbE\int_0^T  \Big[\D l(s) + \lan \D b(s), p(s) \ran_H  + \lan \D \si(s), q(s) \ran_{\cL_2^0} \Big]ds, \\
=& \dbE\int_0^T \Big\{ \lan l_u(t,\Th(t)), v(t) - u^*(t) \ran_U +  \lan l_\mu(t,\Th(t)), v(t-\d) - \mu^*(t) \ran_U \\
&\q + \lan b_u(t,\Th(t))^* p(t),  v(t) - u^*(t) \ran_U +   \lan b_\mu(t,\Th(t))^* p(t), v(t-\d) - \mu^*(t) \ran_U \\
&\q + \lan \si_u(t,\Th(t))^* q(t),  v(t) - u^*(t) \ran_U +   \lan \si_\mu(t,\Th(t))^* q(t), v(t-\d) - \mu^*(t) \ran_U
\Big\} dt.
\end{aligned}$$
Therefore, by the definition of Hamilton function $G(\cd)$,  we derive the  necessary maximum principle condition \rf{mp}.
\ef

\br{}
There are lots of works on the optimal control problem for stochastic delay systems in both finite and infinite dimensional frameworks. We mention a few of them. Chen--Wu \cite{CW} investigated the case of stochastic delay differential equations with pointwise delay in the state and the control. Later Zhang \cite{ZF} extended it by allowing move-average delay in the state and control.
As to the case of stochastic delay evolution equations, Meng--Shen \cite{MS} discussed the problem in the weak solution framework where a version of It\^o's formula is applicable. Under the same framework, Li et al. \cite{LZDXD} studied the stochastic maximum principle in the infinite horizon case. We point out that the final cost in the above papers does not depend on the past of state and the anticipated backward stochastic differential/evolution equations are introduced as the dual equations.
Recently Guatteri et al \cite{GMO} and Guatteri-Masiero \cite{GM} made some progress by allowing the final cost to depend on the past of the state. They introduced certain anticipated backward stochastic integro-differential system due to the dependence on the past trajectory in the final cost.
However, the diffusion term in \cite{GMO} was independent of the state and control.

In contrast with the above papers, we provide a new viewpoint which is inspired by \cite{MSWZ}. The novelties of this method are threefold. First, a new interesting idea via Volterra integral system is introduced for the study of SDEEs. Second, to our best knowledge, it is the first time to show the application of infinite dimensional BSVIEs theory to the optimal control problem for stochastic delayed systems.
Last but not the least, in contrast to previous studies, the final cost functional here is naturally allowed to depend on the past of the state.
\er

\section{Concluding remark}

In this paper, the notion of singular BSVIEs in infinite-dimensional spaces is introduced and the  well-posedness is investigated. We propose a kind of non-convolution singular kernel which can cover many concrete types, and can be well used in treating the forward SVIEs case. To explain the theoretical results, several interesting examples are given in details.  
As applications, we study two classes of optimal control problems. For the case of controlled stochastic delay evolution equation, we use BSVIEs instead of antipated BSDEs as the adjoint equations and derive the maximum principle when the terminal cost depends on the past state for the first time.

\section{Appendix}

In this Appendix, we review two typical papers \cite{LY, Zhang} and compare our results on the forward SVIEs with theirs.

\subsection{Comparison with Lin-Yong \cite{LY}}

 To begin with, we revisit \cite{LY} where the authors considered the finite dimensional valued Volterra equation of the form
	\bel{Y} y(t) = \eta(t) + \int_0^t \frac{f(t,s,y(s))}{(t-s)^{1-\beta}}ds, \quad a.e.\; t \in [0,T].
	\ee
	Here $\beta \in (0,1)$. For readers' convenience, we list the assumptions and the main result in \cite{LY}.\\% a
	
	(H6.1) Let $f: \D \times \dbR^n \to \dbR^n$ be a measurable map, $y \mapsto f(t,s,y)$ being continuously differentiable. There are nonnegative functions $L_0(\cd), L(\cd)$ with
	$$ L_0(\cdot) \in L^{(\frac{2}{1+2\beta } \vee 1) +} (0,T;\dbR_+), \qquad
	L(\cdot) \in L^{(\frac{1}{\beta} \vee 2)+} (0,T;\dbR_+), $$
	%for some $p \ge 1$ (with the convention that $\frac{1}{0} = \infty $)
	such that  %for some $\beta \in (0,1)$,
	$$ |f(t,s,0)|_{\dbR^n} \les L_0(s), \quad (t,s) \in \D, $$
	$$ |f(t,s,y) - f(t,s,y')|_{\dbR^n}  \les L(s) \cd |y-y'|_{\dbR^n}, \quad (t,s) \in \D,\; y, y' \in \dbR^n,$$
	where $\ds L^{s+}(0, T ;\dbR_+)\=\bigcup_{r>s} L^{r}(0, T;\dbR_+), (1 \les s<\infty).$\\

	In \cite{LY}, the authors showed that when (H6.1) holds and $\b \in (0,1)$, for any $\eta(\cd) \in L^2 (0,T;\dbR^n)$,  \rf{Y} admits a unique solution  $y(\cd) \equiv  y(\cd; \eta(\cd)) \in L^2 (0,T;\dbR^n)$.
	\\
	
	Notice that the equation \rf{Y} can be viewed as a special case of \rf{fSVIE}, with
	%$$
	%
	%
	$$\left\{ \begin{aligned}
		&\f(t) \= \eta(t) + \int_0^t \frac{f(t,s,0)}{(t-s)^{1-\beta}}ds, \\
		&A(t,s,y)\=  \frac{f(t,s,y(s))-f(t,s,0)}{(t-s)^{1-\beta}},\\
		&B(t,s,y) \= 0.
	\end{aligned}\right.$$
	Obviously, $\ds A(t,s,0)=0$ and
	$\ds |A(t,s,x) - A(t,s,y)|_{\dbR^n}  \les \frac{L(s)}{(t-s)^{1-\beta}} |x-y|_{\dbR^n}. $

	{Firstly, we claim that $\f(\cd) \in L^2(0,T;\dbR^n)$.}
	In fact, by the assumption on $f(t,s,0)$ in (H6.1)  and Young's convolution inequality \cite[Corollary 2.2]{LY},
	we choose $p=2$ and $q \ges 1,\; 1\les r <\frac{1}{1-\b}$ satisfying $\frac{3}{2} = \frac{1}{q} + \frac{1}{r}$. Then we have
	$$\(\int_0^T \lt| \int_0^t \frac{f(t,s,0)}{(t-s)^{1-\beta}}ds \rt|_{\dbR^n}^2dt\)^\frac{1}{2} \les \(\int_0^{T} \lt| \int_0^t \frac{L_0(s)}{(t-s)^{1-\b}} ds \rt|_\dbR^2 dt \)^\frac{1}{2}
	\les \( \frac{T^{1-r(1-\b)}}{1-r(1-\b)}\)^\frac{1}{r} \lt\|L_0(\cd)\rt\|_{L^q(0,T;\dbR_+)}.$$

	\no Case 1. When $\ds \b \in (0, \frac{1}{2})$, $\ds L_0(\cd) \in L^{\frac{2}{1+2\b}+}(0,T;\dbR_+)$, we can find $\ds r$ close enough to $\frac{1}{1-\b}$ such that
	$$ \frac{1}{q} = \frac{3}{2} - \frac{1}{r} \uparrow \frac{3}{2} - (1-\b) = \frac{1+2\b}{2},$$
	thus $q\downarrow \frac{2}{1+2\b}$, so that $L_0(\cd) \in L^q(0,T;\dbR_+).$ \\

	\no Case 2. When $\ds \b \in [ \frac{1}{2},1)$, $\ds L_0(\cd) \in L^{1+}(0,T;\dbR_+)$, we can find $\ds r$ close enough to $2$ such that
	$$ \frac{1}{q} = \frac{3}{2} - \frac{1}{r} \uparrow \frac{3}{2} - \frac{1}{2} = 1,$$
	thus $q\downarrow 1$, so that $L_0(\cd) \in L^q(0,T;\dbR_+).$
	
	To sum up, in both cases, it follows from the definition of $\f$ that $\f(\cd) \in L^2(0,T;\dbR^n)$.\\

	Secondly, let $\ds K(t,s)\= \frac{L(s)}{(t-s)^{1-\b}} $. We claim that  $$K(\cd,\cd) \in L^2(\D)  \text{ when } \b \in (\frac{1}{2},1) \text{ and } K(\cd,\cd) \notin L^2(\D) \text{ when } \b \in (0,\frac{1}{2}].$$
	
	\no Observe that
	$\ds L(\cdot) \in  L^{(\frac{1}{\beta} \vee 2)+} (0,T;\dbR_+). $
	In the case $\ds \b \in (0,\frac{1}{2}]$, if we choose $L(s) \equiv 1$, we have
	$$ \int_0^T \int_s^T K(t,s)^2 dtds = \int_0^T \int_s^T \frac{1}{(t-s)^{2(1-\b)}}dtds = \i,$$
	which means that $K(\cd,\cd) \notin L^2(\D)$.
	On the other hand, if $\b \in ( \frac{1}{2},1)$, %similarly, if we use H\"older directly, it's to rough.
	then after basic calculations,
	$$ \int_0^T \int_s^T K(t,s)^2 dtds = \int_0^T \int_s^T \frac{L(s)^2}{(t-s)^{2(1-\b)}}dtds =
	\int_0^T L(s)^2 \int_s^T \frac{1}{(t-s)^{2(1-\b)}}dtds < \i,$$
	which indicates $K(\cd,\cd) \in L^2(\D)$.\\
	
	To sum up the above arguments, when our SVIEs reduce to the finite dimensional deterministic system, we have the followig conclusions. The free term in our paper can cover that in \cite{LY}. As to the Lisphcitz function,  our $K_1(\cd,\cd)$ has more general form. As a tradeoff, we only cover the case $\b\in(\frac{1}{2},1)$.
	
\subsection{Comparison with Zhang \cite{Zhang}}

At first, we review the singularity in \cite{Zhang}.
Let $\kappa: \triangle\to\dbR_+$ be a measurable function.
Assume that %for any $T>0$
$$
\int^\cd_0\kappa(\cd,s)d s\in L^\infty(0,T;\dbR_+)
$$
and
\bel{Limit-epsilon}\ba{ll}
\ns\ds
\limsup_{\e\downarrow 0}\left\|\int^{\cdot+\e}_{\cdot}
\kappa(\cdot+\e,s)d s\right\|_{L^\infty(0,T;\dbR_+)}=0.
\ea\ee

\no Denote all the functions satisfying the above two properties by $\sK_0$. The second condition in the definition of $\sK_0$ is to guarantee the existence of resolvent, from which the author derived a singular Gronwall inequality of Volterra type.
In \cite{Zhang}, the author investigated the following stochastic Volterra integral equation in infinite dimensional spaces:
\bel{PP1}
X(t)=g(t)+\int^t_0 b(t,s,X(s))d s+\int^t_0 \si(t,s,X(s))d W(s),\qq t \in [0,T].
\ee

For readers' convenience, we list the following %global Lipschitz and linear growth
conditions imposed
on the coefficients in \cite{Zhang}, and make some slight adjustments in our framework.\\

\no {\bf (H6.2)} Suppose there exist $\kappa_1, \kappa_2 \in \sK_0$, such that \\
(i) For the free term,
$$
\esssup_{t\in[0,T]}\int^t_0[\kappa_1(t,s)+\kappa_2(t,s)]
\cdot\dbE\|g(s)\|_H^2 d s<+\infty,
$$
(ii) For the drift term, for all $(t,s)\in\D$,
$\omega\in\Omega$ and $x\in H$,
$$
\|b(t,s,\omega,x)\|_H\les \kappa_1(t,s)\cdot(\|x\|_H+1)
$$
and
$$
\|b(t,s,\omega,x)-b(t,s,\omega,y)\|_H\les \kappa_2(t,s)\cdot\|x-y\|_H.
$$

\no (iii) For the diffusion term, for all $(t,s)\in\D$,
$\omega\in\Omega$ and $x,y\in H$,
$$
\|\si(t,s,\omega,x)\|^2_{\cL_2^0}\les \kappa_1(t,s)\cdot(\|x\|_H^2+1)
$$
and
$$
\|\si(t,s,\omega,x)-\si(t,s,\omega,y)\|^2_{\cL_2^0}\les \kappa_2(t,s)\cdot\|x-y\|_H^2.
$$
\\

In \cite{Zhang}, under {\bf (H6.2)}, by Picard iteration and Gronwall inequality, the author proved that  there exists a unique %measurable
adapted process $X(\cd)$ satisfying
\rf{PP1}.\\% and some $C_{T,\kappa_1}>0$. \\
% such that
%\be
%\dbE\|X(t)\|^2_H \les C_{T,\kappa_1}\left[\dbE\|g(t)\|_H^2+\esssup_{t\in[0,T]}\int^t_0\kappa_1(t,s)\cdot\dbE\|g(s)\|_H^2 d s\right]\label{PP2}
%\ee
%for almost all $t\in[0,T]$.
%{\color{red} Is the above useful next?}\\

Next, we compare the assumptions in our paper with that in \cite{Zhang}.

$\bullet$ For the free term,
from {\bf (H6.2)}(ii), for all $(t,s)\in\D$,
$\omega\in\Omega$, we see
$$
\|b(t,s,0)\|_H\les \kappa_1(t,s),\q
\|\si(t,s,0)\|^2_{\cL_2^0}\les \kappa_1(t,s).
$$
\no Since $\kappa_1 \in \sK_0$, it is easy to see that
$$ \dbE \int_0^T \(\int_0^t \|b(t,s,0)\|_H ds\)^2dt < \i, \q \dbE \int_0^T \int_0^t \|\si(t,s,0)\|_{\cL_2^0}^2dsdt<\i. $$
Therefore, the free term in our paper is more weaker than that in \cite{Zhang} by Remark \ref{nonzero}.

$\bullet$ Next we look at the diffusion term. From {\bf (H3.3)} and  {\bf (H6.2) }(iii),
\be \notag \begin{aligned}
	&|B(t,s,x) - B(t,s,y)|^2_{\cL_2^0} \les K_2(t,s)^2 |x-y|^2_H, \\
	&\|\si(t,s,x)-\si(t,s,y)\|^2_{\cL_2^0}\les \kappa_2(t,s)\cdot\|x-y\|_H^2,
\end{aligned} \ee
where $K_2(\cd,\cd)$ and $\kappa_2(\cd,\cd)$ satisfies
$$\ds \esssup_{ s\in [0,T]} \int_s^T K_2(t,s)^2dt < \i, \q
\ds \esssup_{ t\in [0,T]}\int_0^t \kappa_2(t,s)ds <\i.$$
%
%\tb{If we define $K_2(t,s)=0$ and $\kappa_2(t,s)=0$ for $0\les t \les s \les T$.}
{According to Gripenberg \cite[Proposition 2.7, pp.231]{Gripenberg}, $K_2(\cd,\cd)^2$ and $\kappa_2(\cd,\cd)$ are called Volterra kernel of type $L^1$ and $L^\i$ on $[0,T]$, respectively. These two type kernels are two important special cases of Volterra kernels. We believe the reason of using different kernels (at least partially) lies in the different solution spaces. In fact, the solution space in \cite{Zhang} is $L_\dbF^\i(0,T;L^2(\O;H))$ for the convenience of large deviation theory, while in our paper the solution belongs to $L_\dbF^2(0,T;H)$ due to the application to optimal control problem in the sequel.

	Another interesting comparison of $K_2(\cd,\cd)$ and $\kappa_2(\cd,\cd)$ lies in condition 2) of $\sL^2(\D)$ and the above condition (\ref{Limit-epsilon}). In \cite{Zhang}, the condition (\ref{Limit-epsilon}) is to guarantee the existence of resolvent, from which a singular Gronwall inequality of Volterra type can be obtained. In our setting, the condition 2) of $K_2(\cd,\cd)$ is used to keep the successful induction on $[0,T]$. Mathematically speaking, the $\e$ in (\ref{Limit-epsilon}) requires the divided subintervals to be larger than a positive constant uniformly, while it is not needed under our condition 2).

$\bullet$ At last, we discuss the drift term. It seems that the Lipschitz functions of $b$ and $\si$ in \cite{Zhang} have to satisfy certain structure and keep some consistency by their proof. In contrast, no special relationship is required between the Lipschitz functions of the drift and diffusion terms in our paper.

 To sum up the above arguments, both \cite{Zhang} and ours consider the solutions of the forward singular SVIEs in an abstract framework. The free term in our paper can cover that in \cite{Zhang}.  Due to the different solution spaces, the Lispscitz functions of the diffusion term in \cite{Zhang} and ours are two different kernels (Volterra kernels of type $L^1$ and $L^\i$). For the drift term, the Lipschitz function in \cite{Zhang} needs to keep consistency with that in the diffusion term, while in our paper, there is no such restriction.


\begin{thebibliography}{9}
\rm	


\bibitem{Abel} N.~H.~Abel, \it Solution de quelques problemes al'aide d'integrales definies, Werke 1, \sl Mag. Naturvidenskaberne, \rm (1823), 10--12. 

\bibitem{Abi1} E.~Abi Jaber, M.~Larsson, S.~Pulido, \it Affine Volterra processes, \sl  Ann. Appl. Probab., \rm 29 (2019), 3155--3200.

\bibitem{Abi2} E.~Abi Jaber, E.~Miller, H.~Pham, \it  Linear-quadratic control for a class of stochastic Volterra equations: solvability and approximation, \sl  Ann. Appl. Probab., \rm 31 (2021), 2244--2274.

\bibitem{Agram-Djehiche} N. Agram, B. Djehiche, \it On a class of reflected backward stochastic Volterra integral equations and related time-inconsistent optimal stopping problems, \sl Systems Control Lett., \rm 155, (2021), 9 pp.

\bibitem{AO} N.~Agram, B.~Oksendal, \it Malliavin calculus and optimal control of stochastic Volterra equations, \sl J. Optim. Theory Appl., \rm 167 (2015), 1070--1094.


\bibitem{AN} A.~Aman, M.~N'Zi, \it Backward stochastic nonlinear Volterra integral equation with local Lipschitz drift, \sl Probab. Math. Statist., \rm 25 (2005), 105--127.

\bibitem{Y 2011} V.~Anh, W.~Grecksch, J.~Yong, \it Regularity of backward stochastic Volterra integral equations in Hilbert spaces, \sl Stoch. Anal. Appl., \rm 29 (2011), 146--168.

\bibitem{BD}  N.~B\"auerle, S.~Desmettre, \it Portfolio optimization in fractional and rough Heston models, \sl SIAM J. Financial Math., \rm 11 (2020), 240--273.

\bibitem{BRG} P.~Beissner, E.~Rosazza Gianin, \it The term structure of Sharpe ratios and arbitrage-free asset pricing in continuous time, \sl Probab. Uncertain. Quant. Risk, \rm 6 (2021), 23--52.

%\bibitem{BM} M.~Berger, V.~Mizel, \it Volterra equations with It\^o integrals, I and II, \sl J. Integral Equations, \rm 2 (1980) 187--245, 319--337.

\bibitem{CW} L.~Chen, Z.~Wu, \it Maximum principle for the stochastic optimal control problem with delay and application, \sl Automatica J. IFAC, \rm 46 (2010), 1074--1080.

%\bibitem{DBX} X.~Dai, W.~Bu, A.~Xiao, \it Well-posedness and EM approximations for non-Lipschitz stochastic fractional integro-differential equations,  \sl J. Comput. Appl. Math., \rm 356 (2019), 377--390.

%\bibitem{DZL}  H.~Dai, J.~Zhou, H.~Li, \it Infinite horizon stochastic maximum principle for stochastic delay evolution equations in Hilbert spaces, \sl Bull. Malays. Math. Sci. Soc., \rm 44 (2021), 3229--3258.

\bibitem{D} H.~T.~Davis, \it Fractional operations as applied to a class of Volterra integral equations, \sl Amer. J. Math.,  \rm 46 (1924),  95--109.

\bibitem{DU} L.~Decreusefond, A.~S.~\"Ust\"unel,  \it Stochastic analysis of the fractional Brownian motion, \sl Potential Anal., \rm 10 (1999), 177--214.

\bibitem{DJ} J.~Djordjevi\'c, S.~Jankovi\'c,  \it On a class of backward stochastic Volterra integral equations, \sl Appl. Math. Lett., \rm 26 (2013), 1192--1197.

\bibitem{EER}  O.~El Euch, M.~Rosenbaum, \it The characteristic function of rough Heston models, \sl Math. Finance, \rm 29 (2019), 3--38.

%\bibitem{CLP} W.~George Cochran, J.~Lee, J.~Potthoff, \it Stochastic Volterra equations with singular kernels, \sl Stochastic Process. Appl., \rm 56 (1995), 337--349.

\bibitem{Gripenberg} G.~Gripenberg, S.-O.~Londen, O.~Staffans, \sl Volterra integral and functional equations, \rm Cambridge University Press, Cambridge, 1990.

\bibitem{GMO}  G.~Guatteri, F.~Masiero, C.~Orrieri, \it Stochastic maximum principle for SPDEs with delay, \sl Stochastic Process. Appl., \rm 127 (2017), 2396--2427.

\bibitem{GM}  G.~Guatteri, F.~Masiero, \it Stochastic maximum principle for problems with delay with dependence on the past through general measures, \sl Math. Control Relat. Fields, \rm 11 (2021), 829--855.



\bibitem{Hamaguchi1} Y.~Hamaguchi, \it Infinite horizon backward stochastic Volterra integral equations and discounted control problems, \sl ESAIM Control Optim. Calc. Var., \rm 27 (2021), 47 pp.

\bibitem{Hamaguchi} Y.~Hamaguchi, T.~Wang, \it  Linear-quadratic stochastic Volterra controls I: Causal feedback strategies, \rm arXiv:2204.08333.

%\bibitem{H} W.~Han, \it Existence, uniqueness and smoothness results for second-kind Volterra equations with weakly singular kernels, \sl J. Integral Equations Appl., \rm 6 (1994), 365--384.

\bibitem{HO} R.~A.~Handelsman, W.~E.~Olmstead, \it Asymptotic solution to a class of nonlinear Volterra integral equations, \sl SIAM J. Appl. Math., \rm 22(1972), 373--384.

\bibitem{Hernandez 2021} C. Hern\'{a}ndez, \it On quadratic multidimensional type-I BSVIEs, infinite families of BSDEs and their applications, \sl Stochastic Process. Appl., \rm 162 (2023), 249--298. 

\bibitem{Hernandez-Possamai 2021} C.~Hern\'andez, D.~ Possamai, \it A unified approach to well-posedness of type-I backward stochastic Volterra integral equations, \sl Electron. J. Probab., \rm 26 (2021), 35 pp.

\bibitem{HuO} Y.~Hu, B.~Oksendal, \it Linear Volterra backward stochastic integral equations, \sl Stochastic Process. Appl.,
\rm 129 (2019), 626--633.

%\bibitem{I}  I.~It\^o, \it On the existence and uniqueness of solutions of stochastic integral equations of the Volterra type, \sl Kodai Math. J., \rm 2 (1979), 158--170.

\bibitem{LP} Q.~Lei, C.~Pun, \it Nonlocal fully nonlinear parabolic differential equations arising in time-inconsistent problems,
\sl J. Differential Equations, \rm 358 (2023), 339--385.


\bibitem{LZDXD}  H.~Li, J.~Zhou, H.~Dai, B.~Xu, W.~Dong,  \it Infinite horizon stochastic delay evolution equations in Hilbert spaces and stochastic maximum principle, \sl Taiwanese J. Math., \rm 26 (2022), 635--665. 

\bibitem{LHH} M.~Li, C.~Huang, Y.~Hu, \it Numerical methods for stochastic Volterra integral equations with weakly singular kernels, \sl IMA J. Numer. Anal., \rm 42 (2022), 2656--2683.

\bibitem{Lin} J. Lin, \it Adapted solution of backward stochastic nonlinear Volterra integral equation, \sl Stoch. Anal. Appl., \rm 20 (2002), 165--183.

\bibitem{LY}  P.~Lin, J.~Yong, \it Controlled singular Volterra integral equations and Pontryagin maximum principle, \sl SIAM J. Control Optim., \rm 58 (2020), 136--164.

\bibitem{L} P.~Linz, \it Numerial methods for Volterra integral equations with singular kernels, \sl SIAM J. Numer. Anal., \rm 6 (1969), 365--374.

%\bibitem{LCE}  C.~E.~Love, \it Singular integral equations of the Volterra type, \sl Trans. Amer. Math. Soc., \rm 15 (1914), 467--476.

\bibitem{LZ} Q.~L\"u, X.~Zhang, \sl  Mathematical control theory for stochastic partial differential equations, \rm Springer, Switzerland, 2021. 

\bibitem{NA}  N.~I.~Mahmudov, A.~Ahmadova, \it Some results on backward stochastic differential equations of fractional order, \sl Qual. Theory Dyn. Syst., \rm 21 (2022), 129 pp.

\bibitem{MS}  Q.~Meng, Y.~Shen, \it Optimal control for stochastic delay evolution equations, \sl Appl. Math. Optim., \rm 74 (2016), 53--89.

%\bibitem{MS2021} W.~Meng, J.~Shi, \it A global maximum principle for stochastic optimal control problems with delay and applications, \sl Systems Control Lett., \rm 150 (2021).

\bibitem{MSWZ} W.~Meng, J.~Shi, T.~Wang, J.~Zhang, \it  A general maximum principle for optimal control of stochastic differential delay systems, \rm submitted to SIAM J. Control Optim., arXiv:2302.03339, 2023.

\bibitem{MF} R.~K.~Miller, A.~Feldstein, \it Smoothness of solutions of Volterra integral equations with weakly singular kernels, \sl SIAM J. Math. Anal., \rm 2 (1971), 242--258.



\bibitem{OR} L.~Overbeck, J.~R\"oder, \it Path-dependent backward stochastic Volterra integral equations with jumps, differentiability
and duality principle, \sl  Probab. Uncertain. Quant. Risk, \rm 3 (2018), 4 pp.


\bibitem{PP} E.~Pardoux, S.~Peng,  \it Adapted solution of backward stochastic differential equation, \sl Systems Control Lett., \rm 14 (1990), 55--61.

\bibitem{Popier 2021} A.~Popier, \it Backward stochastic Volterra integral equations with jumps in a general filtration, \sl ESAIM Probab. Stat., \rm 25 (2021), 133--203.

%\bibitem{PS}  D.~Pr\"omel, D.~Scheffels, \it Stochastic Volterra equations with H\"older diffusion coefficients, \sl Stochastic Process. Appl., \rm 161 (2023), 291--315.

\bibitem{PS}  D.~Pr\"omel, D.~Scheffels, \it Pathwise uniqueness for singular stochastic Volterra equations with H\" older coefficients, \rm  arXiv:2212.08029. 

\bibitem{P} J.~Pr\"uss, \it Evolutionary integral equations and applications, \sl Monogr. Math., \rm Birkh\"auser, 1993.

\bibitem{Ren} Y.~Ren, \it  On solutions of backward stochastic Volterra integral equations with jumps in Hilbert spaces, \sl J. Optim. Theory Appl., \rm 144 (2010), 319--333.

\bibitem{RM}  J.~H.~Roberts, W.~R.~Mann, \it On a certain nonlinear integral equation of the Volterra type, \sl Pacific J. Math., \rm 1 (1951), 431--445.

\bibitem{W 2012} Y.~Shi, T.~Wang, \it  Solvability of general backward stochastic Volterra integral equations, \sl J. Korean Math. Soc., \rm 49 (2012), 1301--1321.

\bibitem{Shi-Wang-Yong 2013} Y.~Shi, T.~Wang, J.~Yong,  \it Mean-field backward stochastic Volterra integral equations, \sl Discrete Contin. Dyn. Syst. Ser. B, \rm 18 (2013), 1929--1967.

\bibitem{SWY} Y.~Shi, T.~Wang, J.~Yong, \it Optimal control problems of forward-backward stochastic Volterra integral equations, \sl Math. Control Relat. Fields, \rm 5 (2015), 613--649.

\bibitem{SWX} Y.~Shi, J.~Wen,  J.~Xiong, \it Backward doubly stochastic Volterra integral equations and their applications, \sl J. Differential Equations, \rm 269 (2020), 6492--6528.

%\bibitem{VZ} F.~Viens, J.~Zhang, \it A martingale approach for fractional Brownian motions and related path dependent PDEs, \sl Ann. Appl. Probab., \rm 29 (2019), 3489--3540.

\bibitem{WH} H.~Wang, \it Extended backward stochastic Volterra integral equations, quasilinear parabolic equations, \sl Stoch. Dyn., \rm 21 (2020), 37 pp.

\bibitem{WSY} H.~Wang, J.~Sun, J.~Yong, \it Recursive utility processes, dynamic risk measures and quadratic
backward stochastic Volterra integral equations, \sl Appl. Math. Optim., \rm 84 (2021), 145--190.

\bibitem{WHY} H.~Wang, J.~Yong, \it Time-inconsistent stochastic optimal control problems and backward stochastic Volterra integral equations,
\it ESAIM Control Optim. Calc. Var., \rm 27 (2021),  40 pp.

\bibitem{WHYZhang} H.~Wang, J.~Yong, J.~Zhang, \it Path dependent Feynman-Kac formula for forward backward stochastic Volterra integral equations,
\sl Ann. Inst. Henri Poincar\'e Probab. Stat., \rm 58 (2022), 603--638.

\bibitem{WHYZhou} H.~Wang, J.~Yong, C.~Zhou, \it Linear-quadratic optimal controls for stochastic Volterra integral equations: causal state feedback and path-dependent Riccati equations,
\sl SIAM J. Control Optim., \rm 61 (2023), 2595--2629.




\bibitem{W} T.~Wang, \it Linear quadratic control problems of stochastic Volterra integral equations, \sl ESAIM Control Optim. Calc. Var., \rm 24 (2018), 1849--1879.

\bibitem{Wang-2022} T. Wang, \it Backward stochastic Volterra integro-differential equations and applications in optimal control problems, \sl SIAM J. Control Optim., \rm 60 (2022),  2393--2419.

\bibitem{WY} T.~Wang, J.~Yong, \it Backward stochastic Volterra integral equation--representation of adapted solutions, \sl Stochastic Process. Appl., \rm 129 (2019), 4926--4964.

\bibitem{Wang-Yong-2023} T.~Wang, J.~Yong, \it Spike variations for stochastic Volterra integral equations, \sl to appear in SIAM J. Control Optim., \rm  2023, Doi: 10.1137/22M1522097.

\bibitem{WZ} T.~Wang, H.~Zhang, \it Optimal control problems of forward-backward stochastic Volterra integral equations with closed control regions, \sl SIAM J. Control Optim., \rm 55 (2017), 2574--2602.

\bibitem{Wang-Zheng-2021} T. Wang, H. Zheng, \it Closed-loop equilibrium strategies for general time inconsistent optimal control problems, \sl SIAM J. Control Optim. \rm 59 (2021), 3152-3178.
%\bibitem{W} Z.~Wang, \it Existence and uniqueness of solutions to stochastic Volterra equations with singular kernels and non-Lipschitz coefficients, \sl Statist. Probab. Lett., \rm 78 (2008), 1062-1071.

\bibitem{Wang-Zhang 2007} Z.~Wang, X.~Zhang, \it  Non-Lipschitz backward stochastic Volterra type equations with jumps, \sl Stoch. Dyn., \rm 7 (2007), 479--496.
%
%\bibitem{WHu} H.~Wu, J.~Hu, \it Mean field backward doubly stochastic Volterra integral equations and their applications, \sl Discrete Contin. Dyn. Syst. Ser. S, \rm 16(2023), 937--966.

\bibitem{YY} W.~Yan, J.~Yong, \it Time-inconsistent optimal control problems and related issues,  \sl Vol. 164 of Modeling, Stochastic Control, Optimization, and Applications, IMA Volumes in Mathematics and Its Applications, Springer, \rm (2019) 533--569.

\bibitem{Yong 06} J.~Yong, \it Backward stochastic Volterra integral equations and some related problems, \sl Stochastic Process. Appl., \rm 116 (2006), 779--795.

\bibitem{Yong 2007} J.~Yong, \it Continuous-time dynamic risk measure by backward stochastic Volterra integral equations,  \sl Appl. Anal., \rm 86 (2007), 1429--1442.

\bibitem{Y 2008} J.~Yong, \it Well-posedness and regularity of backward stochastic Volterra integral equations, \sl Probab. Theory Relat. Fields, \rm 142 (2008), 21--77.

\bibitem{ZF} F.~Zhang, \it Stochastic maximum principle for optimal control problems involving delayed systems, \sl Sci. China Inf. Sci., \rm 64 (2021), 3 pp.

%\bibitem{Zhang2} X.~Zhang, \it Euler schemes and large deviations for stochastic Volterra equations with singular kernel, \sl  J. Differential Equations, \rm 244 (2008), 2226--2250.

\bibitem{Zhang}  X.~Zhang, \it Stochastic Volterra equations in Banach spaces and stochastic partial differential equation, \sl J. Funct. Anal., \rm 258 (2010), 1361--1425.

\bibitem{ZL} J.~Zhou, B.~Liu,  \it The existence and uniqueness of the solution for nonlinear Kolmogorov equations, \sl J. Differential Equations, \rm 253 (2012), 2873--2915.

\bibitem{ZY} Y.~Zhou, J.~Wang, L.~Zhang, \sl Basic theory of fractional differential equations, \rm World scientific, 2016.









































\end{thebibliography}
\end{document}